\let\latexarabic\arabic
\let\latexdocument\document
\let\latexenddocument\enddocument

\documentclass{biometrika}

\let\document\latexdocument
\let\enddocument\latexenddocument
\AtEndDocument{\printhistory} 
\let\arabic\latexarabic 
\def\rm{}  

\usepackage{amsmath}
\usepackage{amssymb}
\usepackage{times}
\usepackage{bm}
\usepackage{natbib}

\usepackage[plain,noend]{algorithm2e}

\makeatletter
\renewcommand{\algocf@captiontext}[2]{#1\algocf@typo. \AlCapFnt{}#2} 
\def\@algocf@capt@plain{top}
\renewcommand{\algocf@makecaption}[2]{%
  \addtolength{\hsize}{\algomargin}%
  \sbox\@tempboxa{\algocf@captiontext{#1}{#2}}%
  \ifdim\wd\@tempboxa >\hsize
    \hskip .5\algomargin%
    \parbox[t]{\hsize}{\algocf@captiontext{#1}{#2}}
  \else%
    \global\@minipagefalse%
    \hbox to\hsize{\box\@tempboxa}
  \fi%
  \addtolength{\hsize}{-\algomargin}%
}
\makeatother

\SetKwBlock{Repeat}{Repeat}{}

\newcommand*\diff{\mathop{}\!\mathrm{d}}
\newcommand{\norm}[1]{\left\lVert#1\right\rVert}

\usepackage{subcaption}
\usepackage{graphicx}
\usepackage{url}

\begin{document}

\jname{Biometrika}
\jyear{2017}
\jvol{103}
\jnum{1}
\accessdate{Advance Access publication on 31 July 2016}

\received{2 January 2017}
\revised{1 April 2017}

\markboth{O. L. Sandqvist}{Landmarking supermodels and boosted trees}

\title{Event history analysis with time-dependent covariates via landmarking supermodels and boosted trees}

\author{O. L. SANDQVIST}
\affil{PFA Pension, Sundkrogsgade 4, Copenhagen, Denmark \email{oliver.sandqvist@outlook.dk}}

\maketitle

\begin{abstract}

We propose a nonparametric method for dynamic prediction in event history analysis with high-dimensional, time-dependent covariates. The approach estimates future conditional hazards by combining landmarking supermodels with gradient boosted trees. Unlike joint modeling or Cox landmarking models, the proposed estimator flexibly captures interactions and nonlinear effects without imposing restrictive parametric assumptions or requiring the covariate process to be Markovian. We formulate the approach as a sieve M-estimator and establish weak consistency. Computationally, the problem reduces to a Poisson regression, allowing implementation via standard gradient boosting software. A key theoretical advantage is that the method avoids the temporal inconsistencies that arise in landmark Cox models. Simulation studies demonstrate that the method performs well in a variety of settings, and its practical value is illustrated through an analysis of primary biliary cirrhosis data.

\end{abstract}


\section{Introduction} 

Event data with high-dimensional, time-dependent covariates are becoming increasingly common as modern information systems continue to advance. Such data arise, for example, in medicine, where electronic patient records capture dynamic health information which is used to guide individualized treatment based on disease progression. For instance, a poor prognosis may trigger treatment intensification and vice versa. Another recent example is the use of wearables in insurance, where real-time health information about policyholders promises a way to obtain more precise risk calculations and fairer prices. This paper focuses on such applications, where the aim is to dynamically predict future events based on the available event and covariate information. Several methods have been proposed for hazard estimation in the presence of time-dependent covariates, namely joint modeling~\citep{Tsiatis.Davidian:2004,Rizopoulos:2012}, landmarking~\citep{Anderson.etal:1983,VanHouwelingen:2007}, and super-efficient estimation~\citep{Nielsen:2000,Mammen.Nielsen:2007,Bagkavos.etal:2025}. We henceforth distinguish between \textit{hazards} -- arising as intensities with respect to the full event–covariate filtration -- and \textit{future conditional hazards}, where covariate information is restricted to a past time point. 

In joint modeling, one estimates the joint distribution of the event and covariate processes. A related strand of work exclusively models the event hazard, usually under the additional assumption that the covariate-dependence is Markovian~\citep{Andersen.Gill:1982,Nielsen.Linton:1995,Dabrowska:1997,Nielsen:1999,Lee.etal:2021}. These hazard-only methods include the Cox model with time-dependent covariates, and can be used for inference based on historic data. However, without a model for the future evolution of the covariates, their predictive scope is limited to an infinitesimal time horizon and cannot yield, for example, full survival distributions. 
In landmarking, subjects are followed after a certain time (referred to as the landmark) to ascertain how the future conditional hazard depends on the status of the covariate at the landmark. In the super-efficient approach, it is assumed that the hazard given the covariate does not explicitly depend on time and only depends on the current value of the covariate process, which is further assumed to be Markov. It exploits the law of iterated expectations to represent the future conditional hazard as a conditional expectation of the hazard at a future time. The resulting expression is then estimated using kernel methods. The method’s name derives from a result showing that a certain two-step kernel estimator of the future conditional hazard is $\sqrt{n}$-consistent — a consequence of the strong homogeneity assumptions, which enable pooling across time rather than exclusively relying on local information.

In spite of the many successful applications of these methods, they each come with their own challenges. Joint modeling is intrinsically a harder statistical problem than estimating the future conditional hazards as these are lower-dimensional functions of the joint distribution. To overcome this challenge, one often resorts to restrictive parametric models, which are vulnerable to misspecification~\citep{Ferrer.etal:2019}. Furthermore, the approach leads to harder computational problems as the computation of most estimands involves numerical integration over the covariate distribution. For landmarking, a central critique~\citep{VanHouwelingen:2007,Suresh.etal:2017,Bagkavos.etal:2025} is related to the consistency condition of~\citet{Jewell.Nielsen:1993}, which says that the law of iterated expectations should hold across times in the population-level model. For the landmarking models studied in the literature, including the popular landmark Cox model, this condition is violated or at the very least places strict assumptions on the covariate process and the parameters of the model~\citep{Suresh.etal:2017}. Another possible downside is that landmarking may lead to considerable data reduction when it is implemented via subsampling~\citep{Maltzahn.etal:2021}. While landmarking supermodels aim to avoid this information loss by pooling data across several landmark times, the resulting increase in dataset size can make estimation computationally intensive~\citep{deSwart.etal:2025}. Furthermore, statistical theory for landmarking supermodels only seems to be available for landmark Cox models ~\citep{VanHouwelingen:2007,Nicolaie.etal:2013}. 
For the super-efficient approach, the use of kernel estimators makes it scale poorly in higher dimensions. Moreover, the asymptotic analysis of this approach is technically demanding, so~\citet{Bagkavos.etal:2025} restrict their attention to a one-dimensional covariate process.

In this paper, we use a landmarking supermodel approach to estimate future conditional hazards. Unlike existing methods, we employ boosted trees, whose flexibility and consistency ensure that the criterion of~\citet{Jewell.Nielsen:1993} is satisfied in the large sample limit, unlike in rigid parametric landmark models. Furthermore, they enable effective handling of high-dimensional covariates. Our theoretical contribution extends the results of~\citet{VanHouwelingen:2007} beyond Cox models and establish weak consistency of boosted-tree landmarking within a sieve estimation framework. The methodology accommodates a wide range of landmarking schemes, providing a way to balance efficient use of the available data  with the computational cost. Another advantage of our framework is that it does not rely on the Markov assumption, but rather develops a way to use Markovian information to consistently estimate the associated future conditional hazards.

\newpage 

\section{Proposed method} \label{sec:method} 

\subsection{Model specification}

 For clarity of exposition, we develop the method for survival data and a random uniform landmarking scheme; the extension to multivariate counting processes and more general landmarking schemes is analogous and given in Section A of the supplementary material. Consider $n$ individuals observed in the interval $[0,T]$. For $i=1,\dots,n$, let $N_i$ be a single-jump counting process for the $i$th individual, let $Y_i$ denote the at-risk indicator, and let $W_i$ be a $p$-dimensional covariate process. The observed data is $\{N_i(u),Y_i(u),W_i(u) : u \leq T, i=1,\dots,n\}$. For $0 \leq s \leq t \leq T$ introduce $\mathcal{G}_{s,t} = \sigma\{ N_i(u),Y_i(u),W_i(s) : u \leq t, i=1,\dots,n \}$ which is a filtration as a function of $t$ but only uses Markovian information about the covariate via its value at time $s$. Assume that $(N_1,Y_1,W_1),\dots,(N_n,Y_n,W_n)$ are i.i.d. and that the future conditional hazard $\lambda(t,s,w)$ exists, i.e.,
\begin{align} \label{eq:compensator}
    E\{\diff N_i(t) \mid \mathcal{G}_{s,t-} \} = \lambda\{t,s,W_i(s)\} Y_i(t) \diff t.
\end{align}
Throughout, we use $E$ and $E_n$ to denote expectation with respect to the true and empirical measure, respectively. Generic versions of random variables are written without the subscript $i$, i.e., $(N,Y,W)$. Multi-indices are written without commas, e.g. $x_{kj}$ denotes the component indexed by $k$ and $j$, not the product $kj$. 

\begin{remark}[At-risk process]
    The at-risk process provides information about the initial state and filtering phenomena such as censoring. Imposing independent filtering in the sense of~\citet{Andersen.etal:1993} lets one recover the future conditional hazards of the unfiltered population.
\end{remark}

\begin{remark}[Relation to hazards]
Write $\bar{W}_i(t) = \{ W_i(u) : u \leq t\}$. The future conditional hazard is determined by the hazard $\alpha$ and the distribution of the covariates via the relation
 $$\lambda\{t,s,W_i(s)\} = E[\alpha\{t, \bar{W}_i(t)\} \mid \mathcal{G}_{s,t-}]$$
 per the innovation theorem (Theorem 3.4 in~\citet{Aalen:1978}).  This relation is the starting point for the proposed estimator in~\citet{Bagkavos.etal:2025} but is not used in the construction of the estimator studied here.
\end{remark}

\subsection{Estimation} \label{subsec:estimation}

For a fixed landmark $s$, the setup is equivalent to a standard survival analysis problem with baseline at time $s$ and baseline covariates $W(s)$. In this case, the empirical log-likelihood is $E_n\{L(s,F)\}$ with
\begin{align*}
    L_i(s,F) = \int_s^T F\{t,s,W_i(s)\} \diff N_i(t) - \int_s^T Y_i(t) \exp[F\{t,s,W_i(s)\}] \diff t
\end{align*}
for future conditional log-hazard $F$ belonging to a suitable function space $\mathcal{F}$. Let $S_{iq}$ with $q=1,\dots,Q$, for a fixed integer $Q \in \mathbb{N}$, be random variables independent of the observed data and assume $S_{iq} \sim \textnormal{Unif}(0,T)$. For each $q$, a landmark dataset is obtained by landmarking at time $S_{iq}$ for subject $i$. The landmarking super dataset is then constructed by stacking these $Q$ datasets. Treating each observation as if they were i.i.d. leads to the empirical criterion function $M_n(F)=E_n\{m(F)\}$ with 
$$m_i(F)=Q^{-1} \sum_{q=1}^Q L_i(S_{iq},F).$$ 
The population criterion function is $M(F)=E\{ m(F)\}=T^{-1}E\left\{\int_0^T L(s,F) \diff s \right\}$.  This is analogous to~\citet{VanHouwelingen:2007} but for random landmark times, a general future conditional log-hazard $F$, and using the full log-likelihood rather than Cox's partial likelihood. To ensure that the above expectations are well-defined and to facilitate asymptotic results, we impose Assumption~\ref{assumption:Bounded}. Assuming the existence of an upper bound is mild due to the finite observation window, and the lower bound is likewise mild since if there are any parts of $\mathcal{D}$ where the future conditional hazards are $0$, they are usually known in advance and hence do not have to be estimated.

\begin{assumption}[Uniform boundedness] \label{assumption:Bounded}
    For all $F\in \mathcal{F}$, the future conditional hazard $\exp(F)$ is bounded within the interval $[\Lambda_L,\Lambda_U] \subset (0,\infty)$ on its $(p+2)$-dimensional domain $\mathcal{D} \subseteq \{(t,s,w) : 0\leq s \leq t \leq T, w \in \mathbb{R}^p \}$.  
\end{assumption}

We consider empirical risk minimization of $M$ over a function space $\mathcal{F}$ and assume that the true future conditional log-hazard belongs to this space, that is $\log \lambda  \in \mathcal{F}$. For $\log \lambda \in \mathcal{F}$ to be plausible, it is desirable to allow $\mathcal{F}$ to be infinite-dimensional and non-compact. However, maximizing $M_n$ over $\mathcal{F}$ may then not be well-defined, and even when it is, the maximizer may be difficult to compute and have undesirable large sample properties, confer with~\citet{Chen:2007}. We therefore proceed via sieve M-estimation where the function space is approximated using an increasing sequence of subspaces $\mathcal{F}_{k} \subseteq \mathcal{F}_{k+1} \subseteq \mathcal{F}$. For a given $k$, we let 
$$\mathcal{F}_k=\left\{ F=\sum_{j=1}^{m_k} \nu_{kj} g_j : g_j \in \mathcal{T}(d_{kj},\mathcal{P}_k), \: \log\Lambda_L \leq F \leq \log \Lambda_U  \right\}$$
where $m_k$ is the number of trees, $\nu_{kj}$ is the learning rate, and $\mathcal{T}(d_{kj},\mathcal{P}_k)$ are trees of depth $d_{kj}$ using a finite partition $\mathcal{P}_k$ of $\mathcal{D}$ such that terminal node regions are unions of sets in $\mathcal{P}_k$. We let $\mathcal{P}_k$ consist of disjoint time-covariate hypercubes $B_{k\ell}$ indexed by $\ell=(\ell_1,\dots,\ell_{p+2})$ such that 
$B_{k\ell} = \prod_{h=1}^{p+2} [b_{kh\ell_h},b_{kh(\ell_h+1)})$
for the possible split points $b_{kh\ell_h}$. Furthermore, we assume that $\mathcal{P}_{k+1}$ is a refinement of $\mathcal{P}_k$ in the sense that the split points in the latter are contained in the split points of the former. We further assume that $(m_k, d_{kj},\nu_{kj})$ are chosen such that $\mathcal{F}_k\subseteq\mathcal{F}_{k+1}$, which for example holds if $m_k$ is non-decreasing and $(d_{kj},\nu_{kj})$ does not depend on $k$. Finally, let $d$ be a metric on $\cup_{k=1}^\infty \mathcal{F}_k$ and define $\mathcal{F}$ as the completion of the union under $d$. An explicit expression for $\mathcal{F}$ is given in Proposition~\ref{prop:approx} under certain assumptions about $(m_k,d_{kj},\nu_{kj})$ and $\mathcal{P}_k$.

\begin{remark}[Interval convention]
    We use intervals on the form $[\cdot,\cdot)$ since off-the-shelf software packages for boosting such as XGBoost~\citep{Chen.Guestrin:2016} use a splitting convention of $<$ for the left branch and $\geq$ for the right branch.  
\end{remark}

The hyperparameters $(m_k,d_{kj},\nu_{kj})$ are standard for boosted trees and provides a way to manage model complexity. Introducing $\mathcal{P}_k$ likewise provides complexity control and makes the space $\mathcal{F}_k$ finite-dimensional, which facilitates the use of sieve estimation theory. Furthermore, since $F \in \mathcal{F}_k$ has the form $F(t,s,w)=\sum_{\ell} c_{k\ell} 1\{(t,s,w) \in B_{k\ell} \}$ for suitable real numbers $c_{k \ell}$, it holds for any $F \in \mathcal{F}_k$ that
\begin{align*}
    L_i(s,F) &= \sum_{\ell} c_{k\ell} O_{ik\ell}(s,T]-\exp(c_{k\ell})E_{ik\ell }(s,T],  \\
    O_{ik\ell}(s,T] &= \int_s^T 1[\{t,s,W_i(s)\} \in B_{k\ell}] \diff N_i(t), \\
    E_{ik\ell}(s,T] &= \int_s^T 1[\{t,s,W_i(s)\} \in B_{k\ell}] Y_i(t) \diff t.
\end{align*}
The latter two objects are referred to as occurrences and exposures, respectively. Thus, up to additive constants that do not depend on $F$, the empirical criterion function $M_n$ restricted to $\mathcal{F}_k$ equals the log-likelihood from a model where $O_{ik\ell}(S_{iq},T]$ for fixed $k$ are i.i.d. with $O_{ik\ell}(S_{iq},T]\sim \mathrm{Poisson}\{\exp(c_{k\ell} + \log E_{ik\ell}(S_{iq},T])\}$. This is analogous to the result of~\citet{Lindsey:1995}, where the counting process log-likelihood is linked to Poisson regression. Note that when two different $B_{k\ell}$ have the same coefficient $c_{k\ell}$, one may sum the occurrences and exposures for those two $\ell$ and obtain a sufficient statistic. This may lead to considerable computational and memory gains. This data reduction is particularly effective when covariates are discrete or when the partitioning $\mathcal{P}_k$ is coarse.

Linking $M_n$ to the Poisson distribution implies that its gradient has a well-known expression, circumventing the challenges addressed in~\citet{Lee.etal:2021}. Furthermore, gradient boosting with trees and Poisson log-likelihood loss is implemented in several off-the-shelf software packages such as XGBoost, avoiding the need for custom implementations as in~\citet{Lee.etal:2021,Pakbin.etal:2025}. 

The proposed estimator $\hat{F}_n$ is obtained by gradient boosting of $M_n$ over $\mathcal{F}_n$ followed by truncation to the interval $[\log\Lambda_L,\log \Lambda_U]$. In Section~\ref{sec:largesample}, we prove $\hat{F}_n$ is consistent under suitable regularity conditions when  boosting follows Algorithm 1 in~\citet{Friedman:2001} with a fixed step size.


\section{Large sample properties} \label{sec:largesample} 

This section contains the main result of the paper, which states that the proposed estimator is weakly consistent in a certain $L^1$ norm. Define the measure
$$\mu(B) = T^{-1} E\left[\int_0^T \int_s^T Y(t)1[\{t,s,W(s)\} \in B] \diff t \diff s \right]$$
such that for any integrable $f$ we have
$$\int f \diff \mu = T^{-1} E\left[ \int_0^T \int_s^T Y(t)f\{t,s,W(s)\} \diff t \diff s \right].$$ 
We shall consider $\mathcal{F}$ as a subset of $L^\infty(\mu)$, and since $\mu$ is a finite measure it may further be considered as a subset of $L^p(\mu)$ for any $1 \leq p \leq \infty$. This means that functions that are $\mu$-almost surely equal are identified as the same function. Hence, the consistency shown in Theorem~\ref{theorem:consistency} only identifies the true future conditional hazard $\log \lambda$ up to a $\mu$-null set.

To address the consistency condition of~\citet{Jewell.Nielsen:1993}, we first state an obvious universal approximation theorem for boosted trees.
\begin{proposition}[Universal approximation theorem] \label{prop:approx}
    For the metric $d(F_1,F_2)=\norm{F_1-F_2}_{\mu,1}$ induced by the $L^1(\mu)$ norm and with suitable choices of the hyperparameters $(m_k,d_{kj},\nu_{kj})$ and partitions $\mathcal{P}_k$, it holds that 
    \[
    \mathcal{F} = \Bigl\{ F : 
    F \text{ is $\mu$-measurable and } 
    \log \Lambda_L \leq F \leq \log \Lambda_U \;\mu\text{-almost surely} \Bigr\}.
    \]
\end{proposition}
The proof of Proposition~\ref{prop:approx} is provided in Section C of the Supplementary Material and follows a standard measure-theoretic approximation argument. Consequently, as long as $\lambda$ exists and is bounded within the interval $[\Lambda_L,\Lambda_U]$, the sieves can be constructed so that it lies in the population model. Since the true future conditional hazard trivially satisfies the law of iterated expectations, and since Proposition~\ref{prop:approx} reveals that we have not imposed any structural assumptions (apart from boundedness) that could exclude it, the consistency condition of~\citet{Jewell.Nielsen:1993} poses no difficulties. This is the same situation as in~\citet{Bagkavos.etal:2025}, except that they do not assume $\lambda$ is bounded above, but do impose various forms of differentiability and time-homogeneity.


As mentioned in Section~\ref{sec:method}, unconstrained empirical risk minimization over $\mathcal{F}$ does generally not lead to estimators with desirable statistical properties. Our approach uses sieves to control model complexity, which is similar to the setting of Theorem 4.1 in~\citet{Biau.Cadre:2021}, where controlled growth of the weak learner class combined with a penalized convex loss-function leads to risk-consistency of the boosting estimator. Early stopping, i.e., halting before convergence to the optimizer, provides an alternative way to control model complexity and obtain consistency~\citep{Zhang.Yu:2005}. In our setting, this is governed by the hyper-parameter $m_k$.

To allow for flexibility in the choice of boosting implementation, we make a high-level assumption that boosting approximately maximizes the empirical criterion function. 
\begin{assumption} \label{assumption:maximize}
The estimator $\hat{F}_n$ is an approximate sieve maximizer
    $$M_n(\hat{F}_n) \geq \sup_{F \in \mathcal{F}_n} M_n(F) - o_P(1).$$
\end{assumption}
In section D of the supplementary material, we show that a simple gradient boosting algorithm satisfies Assumption~\ref{assumption:maximize} under weak conditions. The proof is based on Theorem 3.2 of~\citet{Biau.Cadre:2021}, which shows that gradient boosting following~\citet{Friedman:2001} with infinite iterations minimizes the empirical loss under certain conditions.

\begin{theorem}[Consistency] \label{theorem:consistency}
Suppose that Assumptions~\ref{assumption:Bounded} and~\ref{assumption:maximize} hold. Then $d(\hat{F}_n,\log \lambda )=o_P(1)$ where $d$ is the metric $d(F_1,F_2)=\norm{F_1-F_2}_{\mu,1}$ induced by the $L^1(\mu)$ norm.
\end{theorem}

The theorem shows that the boosted-tree landmarking estimator consistently recovers the entire future conditional hazard surface in $L^1(\mu)$, rather than only pointwise predictions at fixed landmark times. The proof of Theorem~\ref{theorem:consistency} is given in Section E of the Supplementary Material. Establishing convergence rates and asymptotic distribution for boosted tree estimators is a known theoretical challenge and is not addressed here; however, we suggest bootstrapping for inference. The extension to data-driven sieve-spaces $\mathcal{F}_k$ via data-adaptive choices of the hyperparameters is also outside the scope of this paper, but we conjecture that the consistency argument could be adapted to hyperparameters that are selected based on sample-splitting or cross-validation.

\section{Numerical study} \label{sec:numerical} 

\subsection{Setup}

We assess the finite-sample performance of the proposed boosted tree landmarking supermodel for predicting survival. Boosting is implemented via the function \texttt{xgb.train} from the \texttt{XGBoost} package in \texttt{R}~\citep{Chen.etal:2025} using Poisson regression with occurrences and exposures as described in Section~\ref{subsec:estimation}; the \texttt{R} implementation~\citep{R:2025} is available on GitHub (\url{https://github.com/oliversandqvist/Web-appendix-landmark-boosting}). The partitioning $\mathcal{P}_n$ is determined by the numerical precision of the simulated data, except for the time variable, which is partitioned into a grid with step size $0.01$. The proposed model is compared with a standard Cox landmarking supermodel fitted using the \texttt{coxph} function in the \texttt{survival} package in R~\citep{Therneau:2024} and two naive models. The naive models estimate the hazard $\alpha$ and then predict survival using $\hat{\alpha}\{t,W(s)\}$ as a biased estimator of the future conditional hazard. This naive prediction approach, also used as a baseline in~\citet{Schoop.etal:2008}, effectively ignores the possibility of future covariate changes. The hazard $\alpha$ is estimated by Cox regression with time-dependent covariates~\citep{Andersen.Gill:1982} and boosted tree Poisson regression based on the partial likelihood~\citep{Andersen.etal:1993}. The latter approach is analogous to~\citet{Lee.etal:2021} except that we use XGBoost for boosting instead of their proposed algorithm.

Hyperparameters for the boosted tree models are obtained using $5$-fold cross-validation combined with a small grid search, with the exception of the high-dimensional setting (Scenario 3) defined below, where three folds are used for $n=100 \, 000$ due to memory limitations. All hyperparameters except the number of trees are selected based on an independent hold-out set with the same values of $n$ and $Q$. The number of trees is selected via cross validation using the training data. The folds are sampled at the subject level rather than at occurrence and exposure level. When $Q=1$, it would also have been meaningful to sample folds at the occurrence and exposure level, but the generalizability that is estimated then measures how well the model can predict the number of occurrences in unseen parts of the domain $\mathcal{D}$ rather than unseen individuals. However, this fails when $Q > 1$, forcing folds to be sampled at the individual level. The hyperparameters are reported in the supplementary material.

 We simulate data for $n$ independent subjects observed continuously on $[0, T]$ with $T=1$. Simulation of jump times is performed using Lewis’ thinning algorithm from~\citet{Ogata:1981}. To evaluate predictive performance, we vary the sample size $n \in \{100, \dots, 100\,000\}$ and the number of landmarks $Q \in \{1, 2, 5, 10\}$. Performance is measured using the root mean squared error (RMSE) between the estimated and true conditional survival probabilities on an independent test set of subjects. The test set is obtained by sampling subjects and uniform landmark times until we obtain $1000$ test subjects where the event or censoring comes after the landmark. In addition to administrative censoring at $T$, subjects have a constant censoring hazard $\lambda_{C}=0.2$, so around $20\%$ of the subjects are censored before administrative censoring. Changes to the covariates similarly occur with a constant hazard of $\lambda_{W}=2$, leading to an average of around two covariate changes per subject. Events occur according to a hazard $\alpha$ to be specified below. As the implied future conditional hazard does not have a closed form expression, we compute the true probability of survival using $100 \,000$ Monte Carlo simulations per test subject with only administrative censoring, i.e., with $\lambda_C=0$. The proposed estimator of survival for the $j$th test subject at the landmark time $s_j$ is the plug-in estimator
\begin{align} \label{eq:surv}
    \hat{S}_j = \exp \left[ -\int_{s_j}^T \hat{F}_n\{t,s_j,W_j(s_j)\} \diff t \right].
\end{align}
The competing models predict survival analogously, replacing $\hat{F}_n$ with their respective future conditional hazard estimators. To test specific properties of the estimator, we consider three scenarios for the hazard $\alpha$ and the distribution of the covariate process.

\textit{Scenario 1 (Linear Markovian)}. This setting has three covariates and a log-linear hazard
$$\alpha\{t,\bar{W}(t)\}=\exp\{\log(0.3)+0.2 t+0.1W_1(t)+0.3W_2(t)+0.3W_3(t)\}.$$ 
We let $W_1(0)$ and $W_2(0)$ be independent Bernoulli variables with parameter $0.5$ while $W_3(0) \sim N(0.5,0.5)$. When the covariates are updated, $W_1$ stays fixed, $W_2$ is overwritten by an independent Bernoulli variable with parameter $0.5$, and an independent $N(0.5,0.25)$ variable is added to the value of $W_3$.

\textit{Scenario 2 (Nonlinear Non-Markovian)}. Let $\tau(t) = \sup \{ 0 \leq s \leq t : \Delta W(s) \neq 0\}$ be the time of the most recent change to the covariates. This setting has three covariates and the hazard
\begin{align*}
    \alpha\{t,\bar{W}(t)\} &= \exp[\log(0.3)+0.3 |\sin\{\pi t  W_2(t)\}|+0.2 \cos\{W_1(t)\} \\
    & \qquad\quad+0.5 \times 1\{W_1(t) = 1, W_3(t) < 0.5\}+0.3 W_3\{\tau(t)-\}^2]
\end{align*}
with the convention $W_3(0-)=0$. The initial value of the covariates and the updates have the same distribution as Scenario 1.

\textit{Scenario 3 (High-dimensional)}. This setting is identical to Scenario 1 except that 47 covariates with no effect on the hazard are included in the estimation. We let $\{W_4(0),\dots,W_{50}(0)\} \sim N(0,\Sigma)$ for $\Sigma = AA^{\textnormal{T}}$ where $\textnormal{T}$ denotes transpose and $A$ is a 47x47 matrix with each entry being an independent standard normal variable. The same $\Sigma$ is used across all subjects. When covariates are updated, an independent $N(0,\Sigma)$ variable is added to the value of $(W_4,\dots,W_{50})$.

\subsection{Results}

\begin{figure}[ht!]
    \centering
\includegraphics[width=1\linewidth]{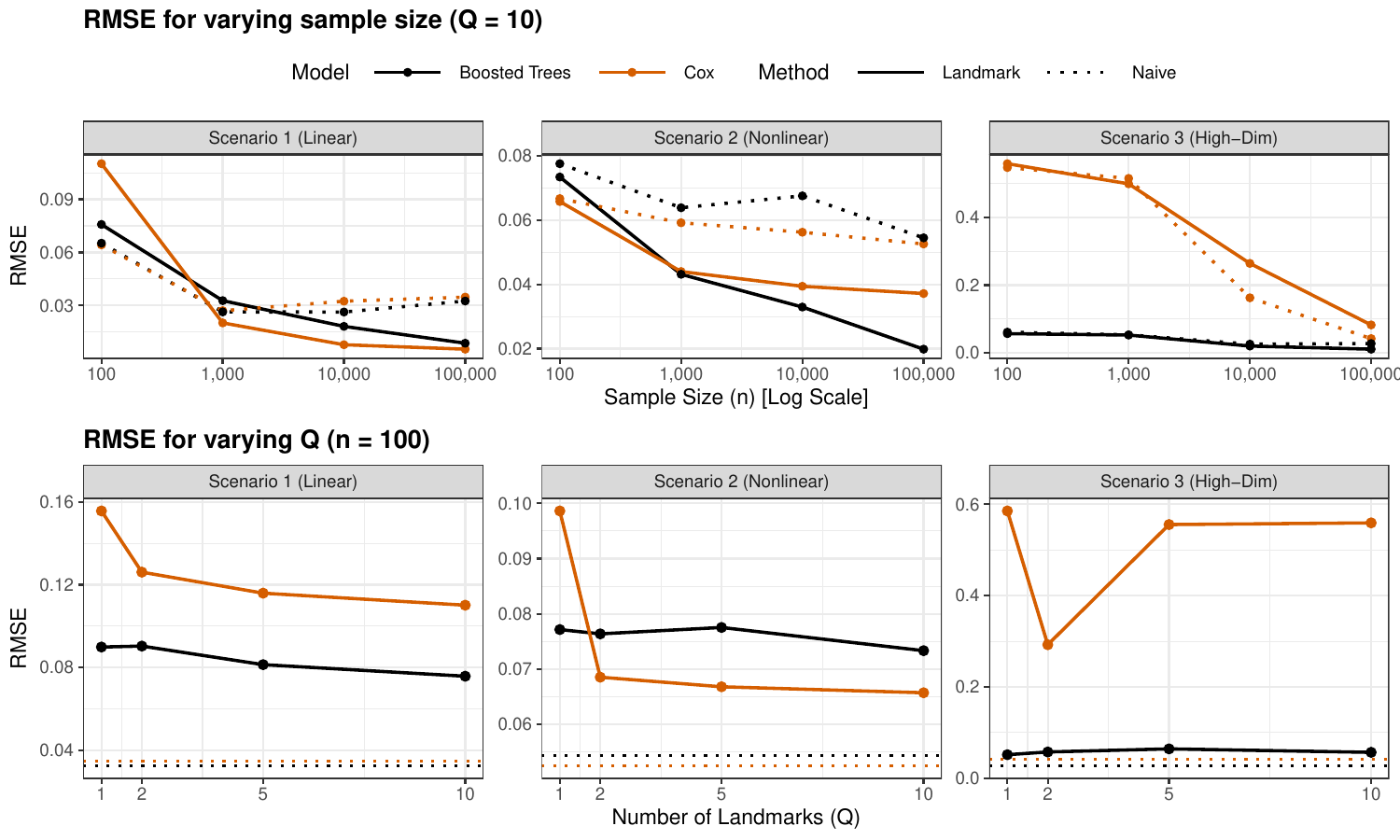}
    \caption{Simulation results. Root mean squared error for the estimated future conditional survival probabilities using the proposed boosted tree landmarking supermodel (black solid lines), the Cox landmarking supermodel (orange solid lines), the naive Cox estimator (orange dotted line), and the naive boosted tree estimator (black dotted line). The top row shows performance as a function of sample size $n$ (log scale) with $Q=10$ landmarks per subject; the bottom row shows performance as a function of $Q$ with sample size $n=100$. In the bottom row, the naive models appear as horizontal lines as they do not landmark the data. The columns correspond to Scenario 1 (Linear Markovian), Scenario 2 (Nonlinear Non-Markovian), and Scenario 3 (High-dimensional).}
    \label{fig:RMSE}
\end{figure}

Fig.~\ref{fig:RMSE} summarizes the results. The top row investigates predictive performance when $Q=10$. In the linear setting (Scenario 1), the naive models perform best for small samples. This is likely because the hazard has a simpler structure than the future conditional hazard, and estimation variance dominates bias at small sample sizes. For large samples, the two landmark supermodels have superior performance and exhibit nearly identical error rates when $n=100 \, 000$. Although the landmark Cox supermodel is misspecified, the deviation from the truth is sufficiently small that it remains the best performing model, as boosted trees require many splits to capture linear effects. In the nonlinear and non-Markovian setting (Scenario 2), the proposed method performs best when $n \geq 1000$, and by a large margin when $n = 100 \, 000$. As expected, the boosted tree estimator significantly outperforms the Cox model in this setting, as the Cox model suffers from asymptotic bias due to its inability to capture the nonlinear structure. In the high-dimensional setting (Scenario 3), the boosted tree estimator demonstrates superior robustness, effectively screening out the 47 noise covariates; in contrast, the unregularized Cox model exhibits high variance at smaller sample sizes. In the supplementary material, we additionally compare the models in terms of mean absolute percentage error (MAPE) with highly similar results. We furthermore compare the models for varying $n$ and fixed $Q \in \{1,2,5\}$, which all show the same pattern as the $Q=10$ case.  

The bottom row of Fig.~\ref{fig:RMSE} examines the impact of the number of landmarks $Q$ for $n=100$. Increasing $Q$ generally leads to improved or comparable performance, likely because the data are used more efficiently. The only exception is the landmark Cox model in Scenario 3, where $Q=2$ is superior. In the supplementary material, we observe that increasing $Q$ continues to result in improved or comparable performance for the boosted trees across the possible values of $n$. The pattern is less clear in the Cox case, where $Q=1$ and $Q=2$ usually perform the best or the worst depending on the scenario and value of $n$. Based on these simulations, choosing $Q=5$ would have provided a good balance between predictive performance and computational cost. The computation time for training the models are reported in the supplementary material.

In summary, the simulations indicate that the proposed method has good predictive performance in a variety of settings, both for small and large samples, as well as for simple and complex data generating mechanisms. The predictive performance improves as the sample size increases and generally improves with an increasing number of landmarks.

\section{Data application} 

The proposed method is demonstrated by an application to the open-source dataset pbc2 from the Mayo Clinic~\citep{Therneau.Grambsch:2000,Fleming.Harrington:1991} to predict clinical progression (i.e., transplantation or death) of patients diagnosed with primary biliary cirrhosis (PBC). The code is available at the GitHub repository listed in Section~\ref{sec:numerical}. The pbc2 dataset contains longitudinal data from a  clinical trial which sought to compare the effectiveness of D-penicillamine versus placebo for treating PBC. It includes patient characteristics (age, sex, and treatment group) and various clinical measurements, from $n=312$ patients. Patients were enrolled during 1974-1984 and observed at subject-specific, irregularly spaced visit times until the minimum date of transplantation, death, or censoring, with administrative censoring occurring at April 30, 1988. At the time of study exit, patients’ statuses were recorded as 143 alive, 29 transplanted, and 140 dead. 

The methodology in this paper is developed for a setting where subjects and their covariates are under continuous observation. In the pbc2 dataset, covariate measurements are only available at visit times. Thus, unless one interpolates the value of the covariates, the only possible landmark times are at these visit times or after study exit where a landmark would contribute nothing to the occurrences and exposures. It was empirically observed (see Fig. 4 in the supplementary material) that the post-enrollment visit times are approximately uniformly distributed between the time of enrollment and the time of study exit. We may hence consider a uniformly distributed landmark time as arising from the following procedure for the $q$'th landmark of the $i$'th subject. Simulate a uniform random variable on $[0,T]$, where $T=14.31$ is the end of the observation window. Time is measured in years since study entry as this is what is registered in the dataset. If this uniform variable exceeds the time of study exit, let $S_{iq}$ be equal to the uniform variable, in which case the landmark contributes nothing to the occurrences and exposures. Alternatively, if the uniform variable is smaller than the time of study exit, simulate $S_{iq}$ from a discrete uniform distribution with atoms at subject $i$'s post-enrollment visit times.

\begin{remark}[Non-uniform landmark distribution]
    If the visit times had not been approximately uniformly distributed, it would be necessary to consider more general distributions with support on $[0,T]$ and perhaps also allow the distribution to depend on the time of study exit; see Section A of the supplementary material for a discussion of how the results of the paper may be generalized to such a setting. We emphasize that when subjects are under continuous observation, one may always simulate $S_{iq}$ uniformly. The use of other distributions may be of interest if they can improve the statistical efficiency of the procedure, but it is only in the case with intermittent observations that this generalization may be necessary in order for the results of the paper to be applicable.  
\end{remark}

We let the partitioning $\mathcal{P}_n$ be determined by the numerical precision of the observations, except for the time variable which is partitioned into a monthly grid. Gradient boosting is performed via the XGBoost package~\citep{Chen.Guestrin:2016} in R~\citep{R:2025}. Hyperparameters are chosen via 5-fold cross-validation with sampling taking place at the subject level rather than at the occurrence and exposure level. Six out of the 17 covariates have at least one missing value (see Table 1 of the supplementary material), but no observations are removed as XGBoost supports missing values by default. There is on average $5.23$ post-enrollment visit times per subject and the average study exit time is $6.41$, so we use $Q=10$ landmark times per subject to increase the chance of each visit time being used at least once. The dataset is not too large, so the estimation runs in about a minute on a standard laptop. 

The native XGBoost importance measure and Shapley values both indicate that the level of  bilirubin in mg/dl (normal range of 0.2 to 1.2) and albumin in g/dl (normal range of 3.4 to 5.4) are strong predictors of clinical progression, see Figs. 5 and 6 of the supplementary material. Figs.~\ref{fig:pdp} and~\ref{fig:marginal} display model predictions for these covariates using partial dependence plots and marginal plots. The partial dependence plots show the fitted mean prediction when the covariate is fixed a given value while averaging over the empirical distribution of the remaining covariates, thereby isolating the fitted effect. The marginal plots group the observed covariate values into $10$ equally spaced bins and depict the average fitted prediction within each bin, thereby illustrating how the predictions vary across the empirical covariate distribution. Averages are computed at the occurrence and exposure level and each graph is supplemented by a rug plot which has a tick at each observation at the occurrence and exposure level. In both plots, the predicted future conditional hazard of clinical progression seems to increase monotonically as a function of bilirubin until around 15 mg/dl where it levels off. For albumin, the plots indicate that the prediction is constant until around 2.5 g/dl after which it decreases monotonically. In both cases, there is a strong dependence outside of the normal range. The marginal plots become irregular for high levels of bilirubin and low levels of albumin since there are few observations there.

\begin{figure}[ht!]
\centering

\begin{minipage}{0.49\textwidth}
  \centering
  \includegraphics[width=\textwidth]{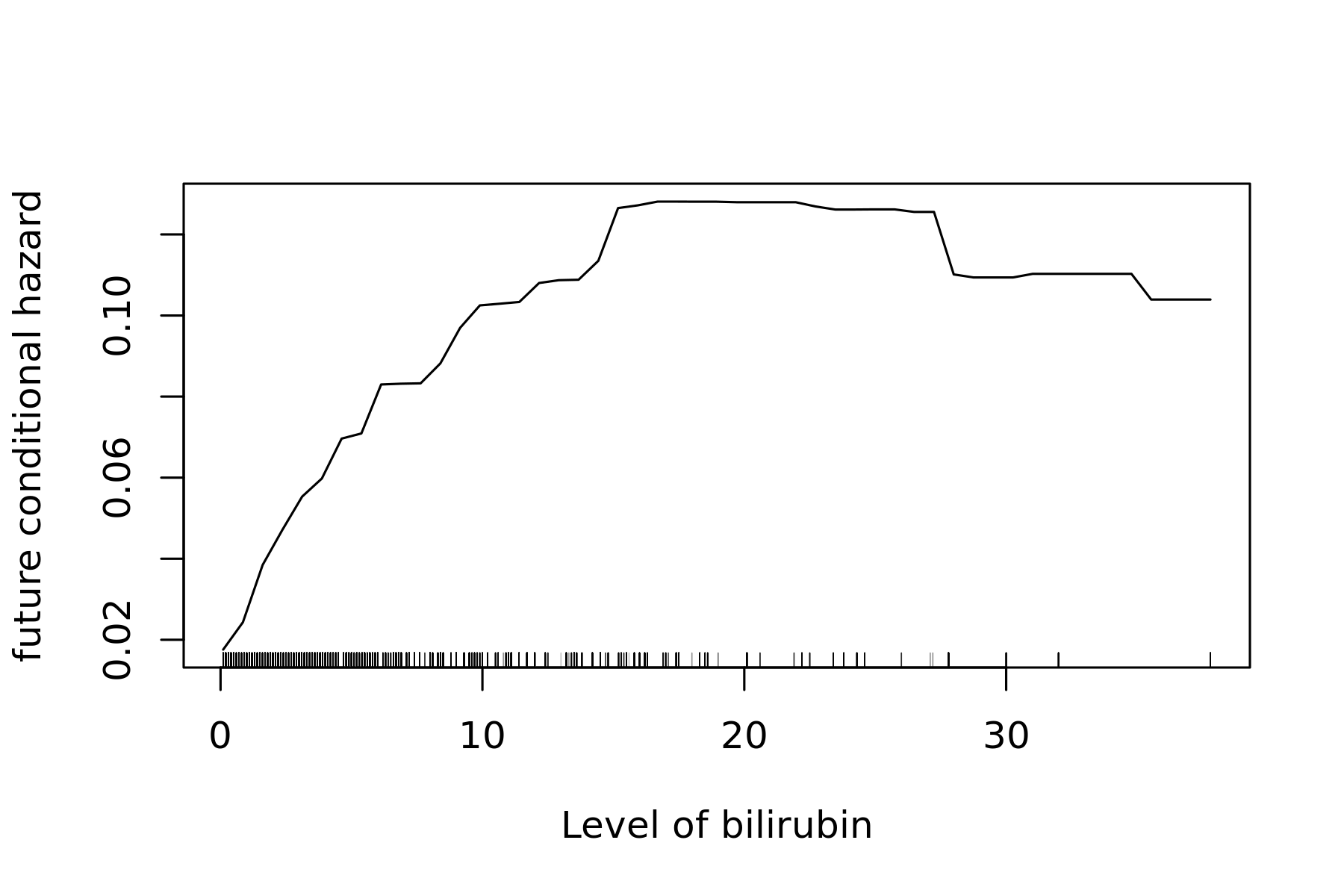}
\end{minipage}
\hfill
\begin{minipage}{0.49\textwidth}
  \centering
  \includegraphics[width=\textwidth]{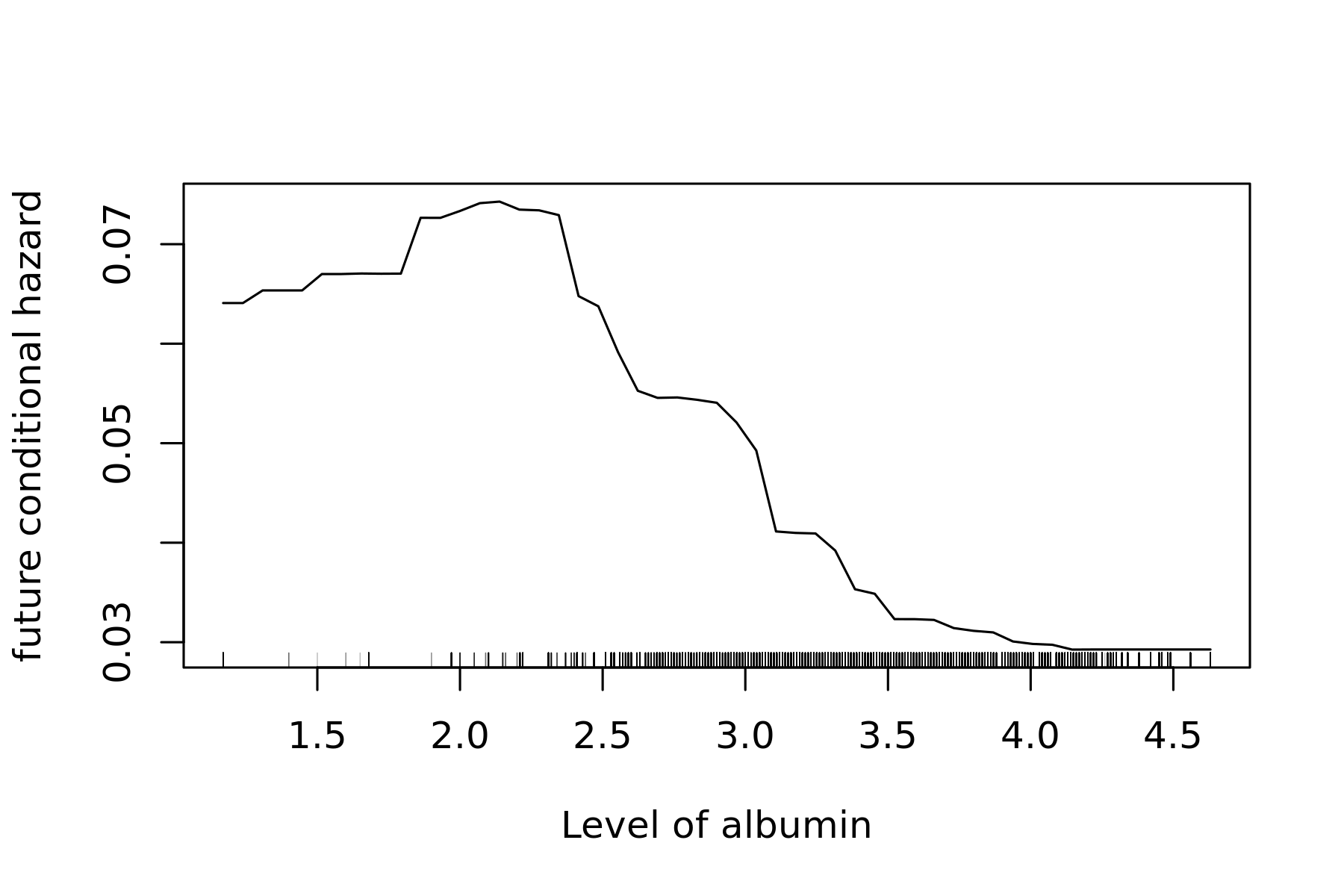}
\end{minipage}
\caption{Partial dependence plot across bilirubin (left) and albumin (right) for the boosted tree estimator (line) and data points (ticks).}
\label{fig:pdp}
\end{figure}

\begin{figure}[ht!]
\centering
\begin{minipage}{0.49\textwidth}
  \centering
  \includegraphics[width=\textwidth]{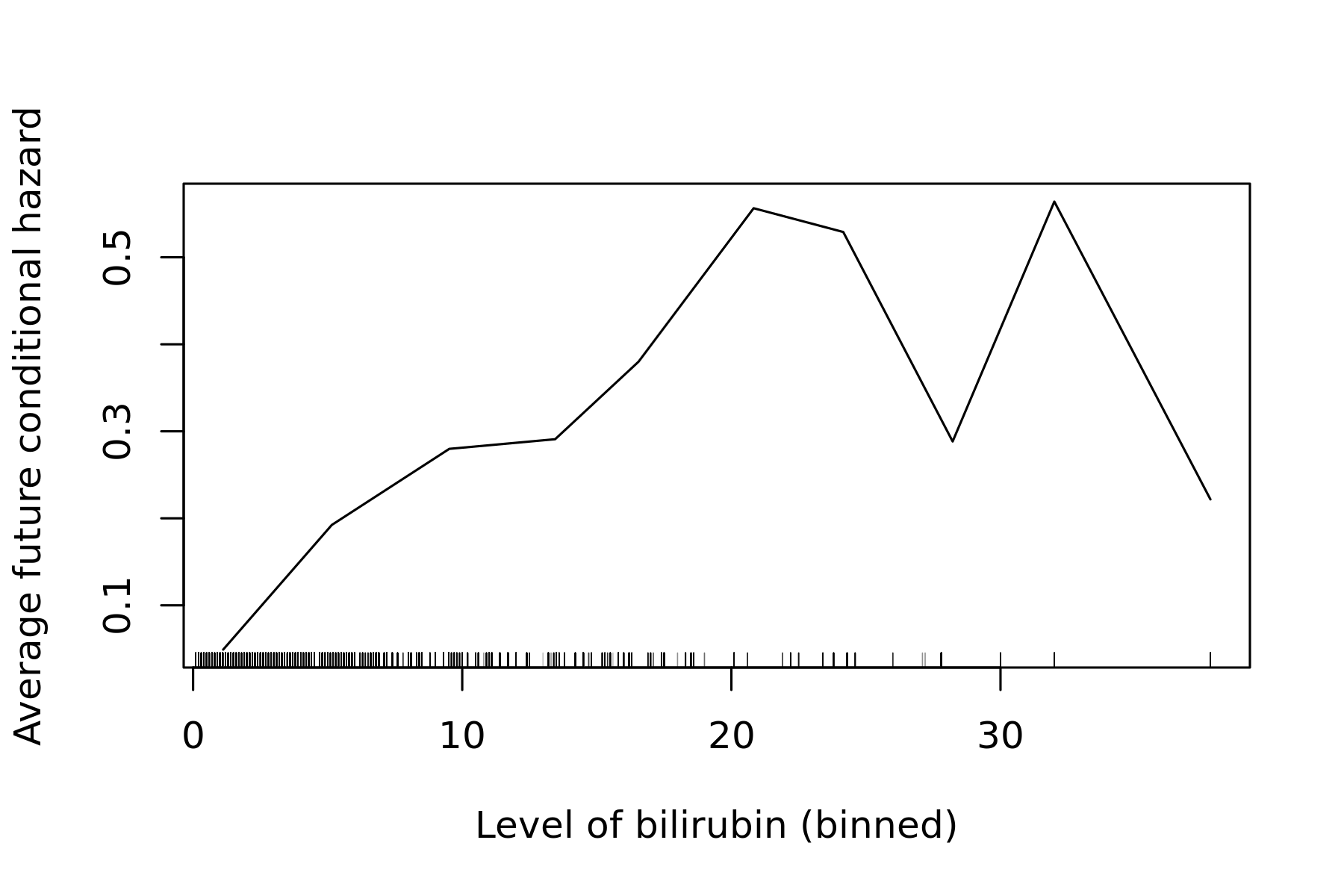}
\end{minipage}
\hfill
\begin{minipage}{0.49\textwidth}
  \centering
  \includegraphics[width=\textwidth]{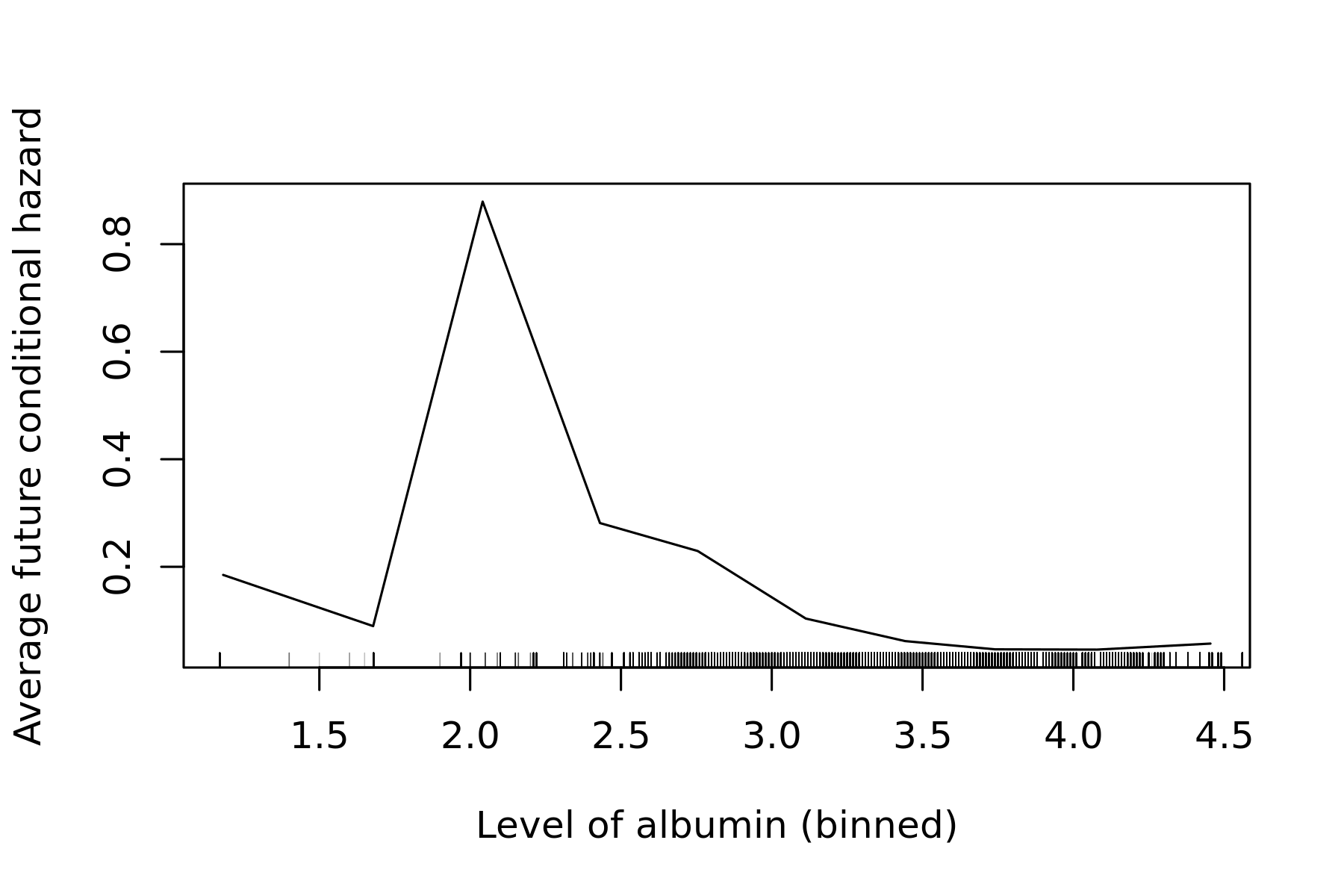}
\end{minipage}
\caption{Marginal plot across bilirubin (left) and albumin (right) for the boosted tree estimator (line) and data points (ticks).}
\label{fig:marginal}
\end{figure}

A main motivation for the methodology is to construct individualized prognoses. We hence compute the subject $j$ conditional survival curve from Equation~\eqref{eq:surv} with $T$ replaced by $t$ for $s_j \leq t \leq s_j+10$ at selected values of $j$ and $s_j$. In Fig. 3, we illustrate the subject-specific survival curve for two representative patients (covariates are listed in Table 2 of the supplementary material). Subject A presented with severe disease (bilirubin 18.5 mg/dL, edema despite diuretics, hepatomegaly) at landmark time $0.85$, whereas Subject B had mild disease (bilirubin 0.5 mg/dL, no edema, no hepatomegaly) at landmark time $0.54$. The contrast between the survival curves is stark, with Subject A and B having predicted year $10$ survival probability of around $0\%$ and $80\%$, respectively. The predicted survival curves are consistent with established results in the literature; see Figs. 1 and 3 of \citet{Dickson.etal:1989}, a widely used and externally validated prognostic model for survival in primary biliary cirrhosis. This agreement supports the ability of the proposed estimator to capture survival differences across individuals and risk strata. 

\begin{figure}[ht!]
\begin{minipage}{0.49\textwidth}
  \centering
\includegraphics[width=\textwidth]{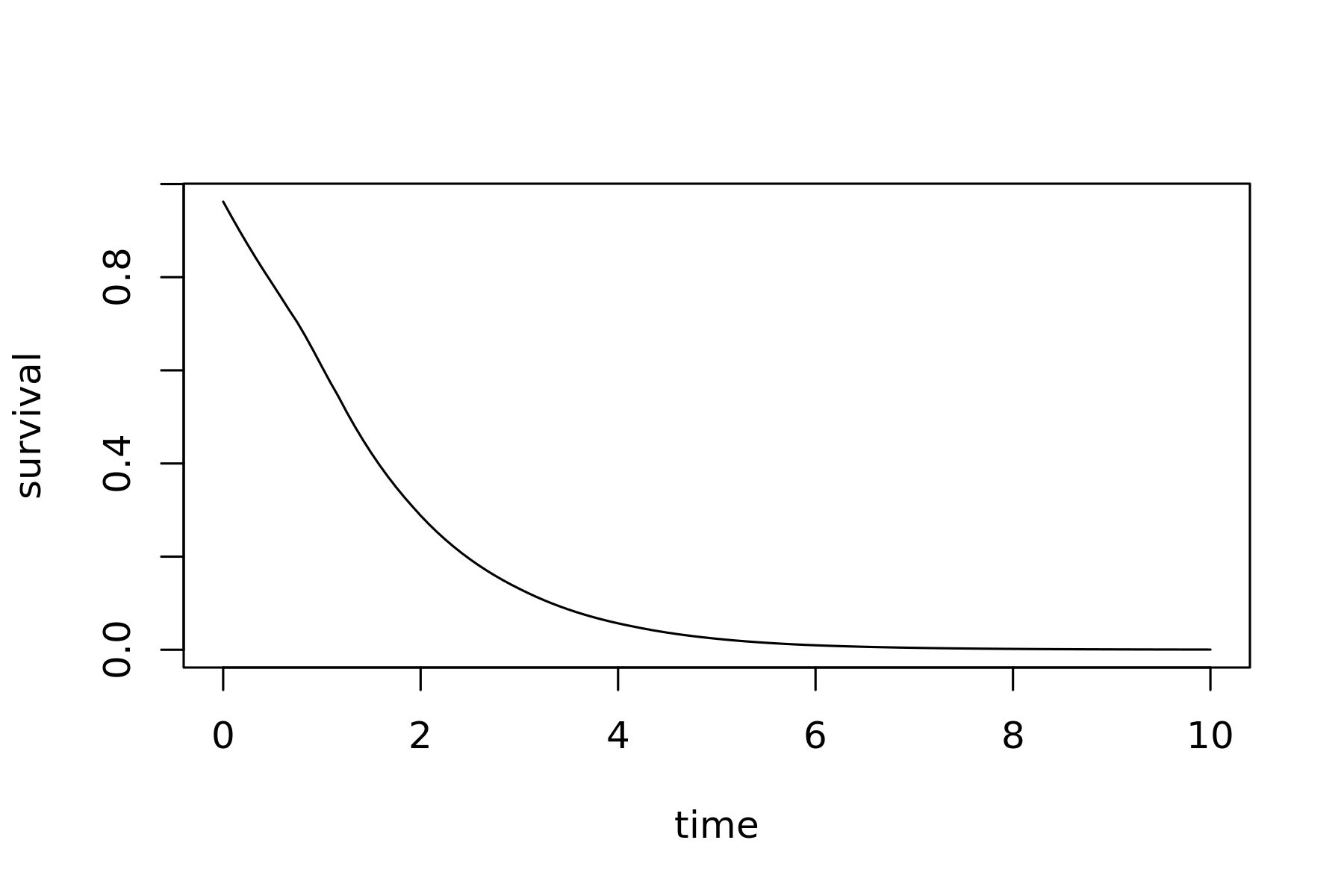}
\end{minipage}
\hfill
\begin{minipage}{0.49\textwidth}
  \centering
\includegraphics[width=\textwidth]{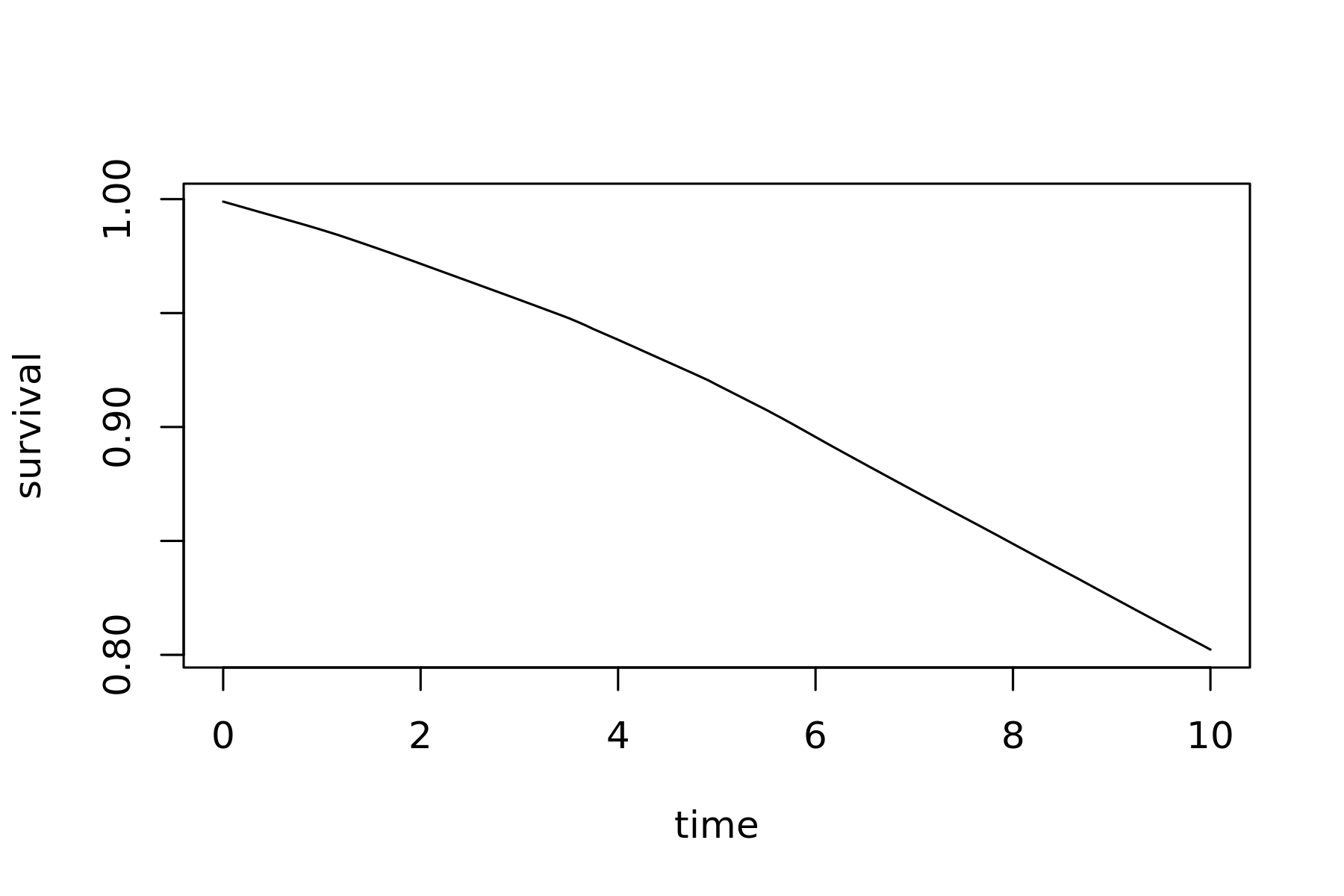}
\end{minipage}
\label{fig:3}
\caption{Survival curve for Subject A (left) and Subject B (right) for the boosted tree estimator.}
\end{figure}

\section{Discussion} 

This paper extends the use of landmarking supermodels to general future conditional hazard models and proposes a boosted tree method for estimation. For ease of exposition, we have restricted our attention to a one-jump counting process with uniform landmark times. This can, however, readily be generalized to multivariate recurrent counting processes and more general landmarking schemes as discussed in Section A of the supplementary material. The simulation study and data application indicate that the method is of practical interest in applied problems where the goal is individualized and dynamical prediction of future events based on the available event and covariate information. 

Our estimated object differs from the kernel-based estimators of~\citet{Putter.Spitoni:2018,Bagkavos.etal:2025,Bladt.Furrer:2025}. The single-landmark approach in~\citet{Putter.Spitoni:2018} is here considered as a kernel-estimator with all its mass placed at a single atom. In the aforementioned papers, the compensator in Equation~\eqref{eq:compensator} is estimated for varying time, but only for a single value of the covariate and the landmark time. In contrast, our approach targets the entire future conditional hazard. This likely yields efficiency gains in integrated error metrics which depend on several values of the covariate and the landmark time, while it may be less efficient at any single landmark and covariate value. It is unclear how these methods could be modified to target the entire future conditional hazard.

The landmark supermodel methodology effectively converts a regression method that is valid for baseline covariates into a regression method that is valid for time-varying covariates. Since the random survival forest proposed in~\cite{Ishwaran.etal:2008,Ishwaran.etal:2014} also estimates a single object for varying baseline covariates, we expect that the landmarking supermodel methodology could be applied to extend their use to time-varying covariates.

An avenue for future research could be to explore whether non-uniform or data-dependent landmarking schemes could increase the statistical efficiency of the approach, and under what conditions. Another avenue might be to generalize the results to penalized criterion functions, which is a salient feature in both~\cite{Chen.Guestrin:2016} and~\cite{Biau.Cadre:2021}. To use standard software implementations like XGBoost, the penalization will have to be applied at the occurrence and exposure level. In order for the approach in this paper to be applicable, one has to ensure that the penalized Poisson objective arises from a criterion function specified at the subject level. 

It might also be interesting to explore the properties of an estimator based on bagging rather than stacking. In the bagging estimator, one would fit a model for each of the $Q$ landmark datasets, and let the final model prediction be the average of the individual predictions. For algorithms whose runtime scales poorly with the size of the dataset, the bagging estimator might be computationally preferable. Furthermore, one would intuitively expect this approach to reduce variance but increase bias compared to the proposed stacking approach.


\section*{Supplementary material}
\label{SM}

The supplementary material includes model extensions, additional details for the simulation and data application, the proofs of Proposition~\ref{prop:approx} and Theorem~\ref{theorem:consistency}, as well as sufficient conditions for Assumption~\ref{assumption:maximize} to hold.


\begin{thebibliography}{7}
\expandafter\ifx\csname natexlab\endcsname\relax\def\natexlab#1{#1}\fi

\bibitem[{Aalen(1978)}]{Aalen:1978}
\textsc{Aalen, O.} (1978).
\newblock Nonparametric inference for a family of counting processes.
\newblock \textit{Ann. Statist.} \textbf{6}, 701--726.

\bibitem[{Andersen et~al.(1993)}]{Andersen.etal:1993}
\textsc{Andersen, P.~K., Borgan, {\O}., Gill, R.~D. \& Keiding, N.} (1993).
\newblock \textit{Statistical models based on counting processes}.
\newblock New York: Springer.

\bibitem[{Andersen \& Gill(1982)}]{Andersen.Gill:1982}
\textsc{Andersen, P.~K. \& Gill, R.~D.} (1982).
\newblock Cox's regression model for counting processes: a large sample study.
\newblock \textit{Ann. Statist.} \textbf{10}, 1100--1120.

\bibitem[{Anderson et~al.(1983)}]{Anderson.etal:1983}
\textsc{Anderson, J.~R., Cain, K.~C. \& Gelber, R.~D.} (1983).
\newblock Analysis of survival by tumor response.
\newblock \textit{J. Clin. Oncol.} \textbf{1}, 710--719.

\bibitem[{Bagkavos et~al.(2025)}]{Bagkavos.etal:2025}
\textsc{Bagkavos, D., Isakson, A., Mammen, E., Nielsen, J.~P. \& Proust–Lima, C.} (2025).
\newblock Super-efficient estimation of future conditional hazards based on time-homogeneous high-quality marker information.
\newblock \textit{Biometrika} \textbf{112}(2), article number asaf008.

\bibitem[{Biau \& Cadre(2021)}]{Biau.Cadre:2021}
\textsc{Biau, G. \& Cadre, B.} (2021).
\newblock Optimization by gradient boosting.
\newblock In \textit{Advances in Contemporary Statistics and Econometrics: Festschrift in Honor of Christine Thomas-Agnan}, pp.~23--44. Cham: Springer.

\bibitem[{Bladt \& Furrer(2025)}]{Bladt.Furrer:2025}
\textsc{Bladt, M. \& Furrer, C.} (2025).
\newblock Conditional Aalen--Johansen estimation.
\newblock \textit{Scand. J. Statist.} \textbf{52}, 873--902.


\bibitem[{Chen(2007)}]{Chen:2007}
\textsc{Chen, X.} (2007).
\newblock Large sample sieve estimation of semi-nonparametric models.
\newblock In \textit{Handbook of Econometrics}, vol.~6, pp.~5549--5632. Amsterdam: Elsevier.

\bibitem[{Chen \& Guestrin(2016)}]{Chen.Guestrin:2016}
\textsc{Chen, T. \& Guestrin, C.} (2016).
\newblock XGBoost: a scalable tree boosting system.
\newblock In \textit{Proc. 22nd ACM SIGKDD Int. Conf. Knowledge Discovery and Data Mining}, pp.~785--794.

\bibitem[{Chen \emph{et al.}(2025)}]{Chen.etal:2025}
\textsc{Chen, T., He, T., Benesty, M., Khotilovich, V., Tang, Y. \& others} (2025).
\newblock \textit{xgboost: Extreme Gradient Boosting}.
\newblock R~package.


\bibitem[{Dabrowska(1997)}]{Dabrowska:1997}
\textsc{Dabrowska, D.~M.} (1997).
\newblock Smoothed Cox regression.
\newblock \textit{Ann. Statist.} \textbf{25}, 1510--1540.

\bibitem[{de~Swart et~al.(2025)}]{deSwart.etal:2025}
\textsc{de~Swart, W.~K., Loog, M. \& Krijthe, J.~H.} (2025).
\newblock A comparative study of methods for dynamic survival analysis.
\newblock \textit{Front. Neurol.} \textbf{16}, article number 1504535.

\bibitem[{Dickson \emph{et al.}(1989)}]{Dickson.etal:1989}
\textsc{Dickson, E.~R., Grambsch, P.~M., Fleming, T.~R., Fisher, L.~D. \& Langworthy, A.} (1989).
\newblock Prognosis in primary biliary cirrhosis: model for decision making.
\newblock \textit{Hepatology} \textbf{10}, 1--7.

\bibitem[{Ferrer et~al.(2019)}]{Ferrer.etal:2019}
\textsc{Ferrer, L., Putter, H. \& Proust--Lima, C.} (2019).
\newblock Individual dynamic predictions using landmarking and joint modelling: validation of estimators and robustness assessment.
\newblock \textit{Stat. Methods Med. Res.} \textbf{28}, 3649--3666.

\bibitem[{Fleming \& Harrington(1991)}]{Fleming.Harrington:1991}
\textsc{Fleming, T.~R. \& Harrington, D.~P.} (1991).
\newblock \textit{Counting Processes and Survival Analysis}.
\newblock New York: Wiley.

\bibitem[{Friedman(2001)}]{Friedman:2001}
\textsc{Friedman, J.~H.} (2001).
\newblock Greedy function approximation: a gradient boosting machine.
\newblock \textit{Ann. Statist.} \textbf{29}, 1189--1232.

\bibitem[{Ishwaran \emph{et al.}(2008)}]{Ishwaran.etal:2008}
\textsc{Ishwaran, H., Kogalur, U.~B., Blackstone, E.~H. \& Lauer, M.~S.} (2008).
\newblock Random survival forests.
\newblock \textit{Ann. Appl. Statist.} \textbf{2}, 841--860.

\bibitem[{Ishwaran \emph{et al.}(2014)}]{Ishwaran.etal:2014}
\textsc{Ishwaran, H., Gerds, T.~A., Kogalur, U.~B., Moore, R.~D., Gange, S.~J. \& Lau, B.~M.} (2014).
\newblock Random survival forests for competing risks.
\newblock \textit{Biostatistics} \textbf{15}, 757--773.


\bibitem[{Jewell \& Nielsen(1993)}]{Jewell.Nielsen:1993}
\textsc{Jewell, N.~P. \& Nielsen, J.~P.} (1993).
\newblock A framework for consistent prediction rules based on markers.
\newblock \textit{Biometrika} \textbf{80}, 153--164.

\bibitem[{Lee et~al.(2021)}]{Lee.etal:2021}
\textsc{Lee, D.~K.~K., Chen, N. \& Ishwaran, H.} (2021).
\newblock Boosted nonparametric hazards with time-dependent covariates.
\newblock \textit{Ann. Statist.} \textbf{49}, 2101.

\bibitem[{Lindsey(1995)}]{Lindsey:1995}
\textsc{Lindsey, J.~K.} (1995).
\newblock Fitting parametric counting processes by using log-linear models.
\newblock \textit{J. R. Statist. Soc. {\rm C}} \textbf{44}, 201--212.

\bibitem[{Maltzahn \emph{et al.}(2021)}]{Maltzahn.etal:2021}
\textsc{Maltzahn, N., Hoff, R., Aalen, O.~O., Mehlum, I.~S., Putter, H. \& Gran, J.~M.} (2021).
\newblock A hybrid landmark Aalen--Johansen estimator for transition probabilities in partially non-Markov multi-state models.
\newblock \textit{Lifetime Data Anal.} \textbf{27}, 737--760.


\bibitem[{Mammen \& Nielsen(2007)}]{Mammen.Nielsen:2007}
\textsc{Mammen, E. \& Nielsen, J.~P.} (2007).
\newblock A general approach to the predictability issue in survival analysis with applications.
\newblock \textit{Biometrika} \textbf{94}, 873--892.

\bibitem[{Nicolaie et~al.(2013)}]{Nicolaie.etal:2013}
\textsc{Nicolaie, M.~A., Van Houwelingen, J.~C., De Witte, T.~M. \& Putter, H.} (2013).
\newblock Dynamic prediction by landmarking in competing risks.
\newblock \textit{Stat. Med.} \textbf{32}, 2031--2047.

\bibitem[{Nielsen \& Linton(1995)}]{Nielsen.Linton:1995}
\textsc{Nielsen, J.~P. \& Linton, O.~B.} (1995).
\newblock Kernel estimation in a nonparametric marker dependent hazard model.
\newblock \textit{Ann. Statist.} \textbf{23}, 1735--1748.

\bibitem[{Nielsen(1999)}]{Nielsen:1999}
\textsc{Nielsen, J.~P.} (1999).
\newblock Super-efficient hazard estimation based on high-quality marker information.
\newblock \textit{Biometrika} \textbf{86}, 227--232.

\bibitem[{Nielsen(2000)}]{Nielsen:2000}
\textsc{Nielsen, J.~P.} (2000).
\newblock Super-efficient prediction based on high-quality marker information.
\newblock \textit{ASTIN Bull.} \textbf{30}, 295--303.

\bibitem[{Ogata(1981)}]{Ogata:1981}
\textsc{Ogata, Y.} (1981).
\newblock On Lewis' simulation method for point processes.
\newblock \textit{IEEE Trans. Inform. Theory} \textbf{27}, 23--31.

\bibitem[{Pakbin et~al.(2025)}]{Pakbin.etal:2025}
\textsc{Pakbin, A., Wang, X., Mortazavi, B.~J. \& Lee, D.~K.~K.} (2025).
\newblock BoXHED2.0: scalable boosting of dynamic survival analysis.
\newblock \textit{J. Stat. Softw.} \textbf{113}, 1--26.

\bibitem[{Putter \& Spitoni(2018)}]{Putter.Spitoni:2018}
\textsc{Putter, H. \& Spitoni, C.} (2018).
\newblock Non-parametric estimation of transition probabilities in non-Markov multi-state models: the landmark Aalen--Johansen estimator.
\newblock \textit{Stat. Methods Med. Res.} \textbf{27}, 2081--2092.

\bibitem[{R Core Team(2025)}]{R:2025}
\textsc{R Core Team} (2025).
\newblock \textit{R: A Language and Environment for Statistical Computing}.
\newblock Vienna: R Foundation for Statistical Computing.

\bibitem[{Rizopoulos(2012)}]{Rizopoulos:2012}
\textsc{Rizopoulos, D.} (2012).
\newblock \textit{Joint models for longitudinal and time-to-event data: with applications in R}.
\newblock Boca Raton, FL: CRC Press.

\bibitem[{Schoop \emph{et al.}(2008)}]{Schoop.etal:2008}
\textsc{Schoop, R., Graf, E. \& Schumacher, M.} (2008).
\newblock Quantifying the predictive performance of prognostic models for censored survival data with time-dependent covariates.
\newblock \textit{Biometrics} \textbf{64}, 603--610.

\bibitem[{Suresh et~al.(2017)}]{Suresh.etal:2017}
\textsc{Suresh, K., Taylor, J.~M.~G., Spratt, D.~E., Daignault, S. \& Tsodikov, A.} (2017).
\newblock Comparison of joint modeling and landmarking for dynamic prediction under an illness-death model.
\newblock \textit{Biometr. J.} \textbf{59}, 1277--1300.

\bibitem[{Therneau(2024)}]{Therneau:2024}
\textsc{Therneau, T.~M.} (2024).
\newblock \textit{A Package for Survival Analysis in R}.
\newblock R~package.

\bibitem[{Therneau \& Grambsch(2000)}]{Therneau.Grambsch:2000}
\textsc{Therneau, T.~M. \& Grambsch, P.~M.} (2000).
\newblock \textit{Modeling Survival Data: Extending the Cox Model}.
\newblock New York: Springer-Verlag.

\bibitem[{Tsiatis \& Davidian(2004)}]{Tsiatis.Davidian:2004}
\textsc{Tsiatis, A.~A. \& Davidian, M.} (2004).
\newblock Joint modeling of longitudinal and time-to-event data: an overview.
\newblock \textit{Statist. Sinica} \textbf{14}, 809--834.

\bibitem[{Van Houwelingen(2007)}]{VanHouwelingen:2007}
\textsc{Van Houwelingen, H.~C.} (2007).
\newblock Dynamic prediction by landmarking in event history analysis.
\newblock \textit{Scand. J. Statist.} \textbf{34}, 70--85.

\bibitem[{Zhang \& Yu(2005)}]{Zhang.Yu:2005}
\textsc{Zhang, T. \& Yu, B.} (2005).
\newblock Boosting with early stopping: convergence and consistency.
\newblock \textit{Ann. Statist.} \textbf{33}, 1538--1579.

\end{thebibliography}

\end{document}


\markboth{O.L.~Sandqvist}{Landmarking supermodels and boosted trees}

\title{Supplementary material for `Event history analysis with time-dependent covariates via landmarking supermodels and boosted trees'}

\author{O.L.~SANDQVIST}
\affil{PFA Pension, Sundkrogsgade 4, Copenhagen, Denmark 
\email{oliver.sandqvist@outlook.dk}}

\maketitle

The supplementary material is organized as follows. Supplement{\scshape~\ref{supplement:extensions}} discusses extensions of the methodology to multivariate recurrent counting process and more general landmarking schemes. Supplement{\scshape~\ref{supplement:numDetails}} provides additional details for the numerical study. Supplement{\scshape~\ref{supplement:dataDetails}} provides additional details for the data application. Supplement{\scshape~\ref{supplement:space}} contains the proof of Proposition 1. Supplement{\scshape~\ref{supplement:maximize}} presents a boosting algorithm and regularity conditions that ensure Assumption 2 is satisfied. Supplement{\scshape~\ref{supplement:consistency}} contains the proof of Theorem 1.

\section{Model extensions} \label{supplement:extensions}

It is possible to extend the methodology of the paper to multivariate one-jump counting processes $\mathbf{N}_i=(N_{i1},...,N_{iJ})$ and at-risk indicators $\mathbf{Y}_i=(Y_{i1},...,Y_{iJ})$. One relatively straightforward approach is to use the setup from Section~\ref{sec:method} separately for each coordinate $N_{ij}$ to estimate their future conditional hazard $\lambda_j$. Relevant information about jumps of the other counting processes can be included in $W_i$. However, there can be a loss of statistical efficiency arising from estimating each $\lambda_j$ separately instead of jointly if there is some shared structure between them e.g., similar dependence on age across different causes of mortality in a competing risk model. 

One may further extend the methodology of the paper to recurrent counting processes. Let $\bar{N}_{ij}(t) = \{N_{ij}(u) : u \leq t\}$.  The only change to the M-estimation setup in this case is that all $F \in \mathcal{F}$ evaluated at time $t$ now also depend on $\bar{N}_{ij}(t-)$. The proposed tree-based estimation procedure is, however, not always applicable, since standard tree-based methods are designed to take finite-dimensional covariates as inputs and a recurrent counting process may have no finite upper bound on the amount of events that can occur before time $T$. To accommodate this, one could either (1) impose Markov-type assumptions (Markov, semi-Markov, etc.) under which the future conditional hazard only depends on a subset of the jumps of $\bar{N}_{ij}$, (2) learn sufficient representations of the path-information (e.g. Markov) based on data and subsequently apply the tree regression method, or (3) develop gradient boosting methods that are tailored to path-valued covariates. The first approach is simple and is probably a good choice in most practical applications.

To accommodate borrowing structure across estimation of $\lambda_j$ for $j=1,...,J$, one may redefine $\mathcal{G}_{s,t}= \sigma\{ \mathbf{N}_i(u),\mathbf{Y}_i(u),W_i(s) : u \leq t, i=1,\dots,n \}$ which is generated by all the counting and at-risk processes.  The future conditional hazards are then collectively given by
$$E\{\diff \mathbf{N}_{i}(t) \mid \mathcal{G}_{s,t-} \} = \boldsymbol{\lambda}\{t,s,\bar{\mathbf{N}}_i(t-),W_i(s)\} \mathbf{Y}_i(t) \diff t$$
for $\bar{\mathbf{N}}_i=(\bar{N}_{i1},...,\bar{N}_{iJ})$ and $\boldsymbol{\lambda} = (\lambda_1,...,\lambda_J)$, and where the right-hand side is to be read as the Hadamard product (element-wise product). In this case, the log-likelihood becomes 
\begin{align*}
    L_i(s,F) &= \sum_{j=1}^J \int_s^T  F_j\{t,s,\bar{\mathbf{N}}_i(t-),W_i(s)\} \diff N_{ij}(t) \\
    & \qquad \quad- \int_s^T  Y_{ij}(t) \exp[F_j\{t,s,\bar{\mathbf{N}}_i(t-),W_i(s)\}] \diff t
\end{align*}
for $F=(F_1,...,F_J)$. The proofs and software implementations can be carried over to this setting without any notable difficulties.

Note that for the specification of the future conditional hazards and the criterion function, it is not necessary to restrict the $W$-information to Markovian information. One could equally well have conditioned on the path $\bar{W}_i(s)$ instead of $W_i(s)$. This is subject to the same challenges as the path-valued variable $\bar{N}_{ij}$ when used in combination with tree-based estimation.

It is also easy to extend the setup to more general landmarking schemes. Instead of uniformly distributed landmark times, we may let $S_{iq} \sim \nu$ independent of the observed data for any probability measure $\nu$ with support on a subset of $[0,T]$. The only change to the M-estimation setup is that the population criterion function becomes $E[\int_0^T L(s,F) \diff \nu(s)]$. This includes deterministic landmark times as a special case. The theorems and proofs would be identical in this setting.

It is also possible to allow the landmark times to depend on the observed data, but it is outside the scope of this paper to examine the large sample properties of the resulting estimator. An example could be simulating uniformly on $[0,\tau]$ rather than $[0,T]$ for an absorption time $\tau$, e.g., the time of death. The motivation is that $L(s,F)=0$ for $s > \tau$, so when landmark times belong to $(\tau,T]$ they contribute nothing to the occurrences and exposures which seems computationally inefficient. Nevertheless, this leads to a different population criterion function $\tilde{M}(F)=E[\tau^{-1} \int_0^{\tau} L(s,F) \diff s]$ which increases the weight for subjects with short event times compared to the criterion function used in this paper $T^{-1}E[\int_0^\tau L(s,F) \diff s]$ where the contribution is proportional to the survival time. We conjecture that such criterion functions also result in consistent estimators, but that the estimators make different statistical efficiency trade-offs, as the criterion functions encode different notions of optimality.

\section{Additional details for numerical study} \label{supplement:numDetails}

\subsection{Hyperparameters}
The most influential hyperparameters were observed to be the \texttt{eta} (Learning rate), \texttt{max\_depth} (Tree depth), and \texttt{nrounds} (Number of trees). Other hyperparameters used provide small changes to the boosting procedure, namely \texttt{min\_child\_weight} (Minimum sum of hessian needed in a child), \texttt{subsample} (Subsample ratio of the training instances), \texttt{colsample\_bytree} (Subsample ratio of columns when constructing each tree), and \texttt{alpha} (L1 regularization term on weights). For simplicity, we only considered a constant learning rate and tree depth, except for the first tree, where we chose a learning rate of $1$ and a tree depth of $0$. Thus, the first model predicts the overall occurrence-exposure rate, so the boosting algorithm does not have to spend a lot of trees in the beginning with a small learning rate to get predictions on the correct order of magnitude. This results in models with fewer trees and equal predictive performance compared to situations where the learning rate and depth are kept fixed throughout.

In Scenario 1 (Linear Markovian), the hyperparameters used for the naive and landmark boosted tree models are:
\begin{itemize}
    \item \texttt{eta} = 0.1
    \item \texttt{max\_depth} = 1
    \item \texttt{min\_child\_weight} = 20
    \item \texttt{subsample} = 0.9
    \item \texttt{colsample\_bytree} = 0.7
\end{itemize}
In Scenario 2 (Nonlinear Non-Markovian), the hyperparameters used for the naive and landmark boosted tree models are:
\begin{itemize}
    \item \texttt{eta} = 0.1
    \item \texttt{max\_depth} = 3
    \item \texttt{min\_child\_weight} = 20
    \item \texttt{subsample} = 0.9
    \item \texttt{colsample\_bytree} = 0.7
\end{itemize}
In Scenario 3 (High-dimensional), the hyperparameters used for the naive boosted tree models are:
\begin{itemize}
    \item \texttt{eta} = 0.1
    \item \texttt{max\_depth} = 1
    \item \texttt{min\_child\_weight} = 20
    \item \texttt{subsample} = 0.7
    \item \texttt{alpha} = 100
\end{itemize}
and the landmarked boosted tree models used hyperparameters:
\begin{itemize}
    \item \texttt{eta} = 0.1
    \item \texttt{max\_depth} = 1
    \item \texttt{min\_child\_weight} = 100
    \item \texttt{subsample} = 0.7
    \item \texttt{alpha} = 100
\end{itemize}

\subsection{Additional figures}

\begin{figure}[ht!]
    \centering
\includegraphics[width=1\linewidth]{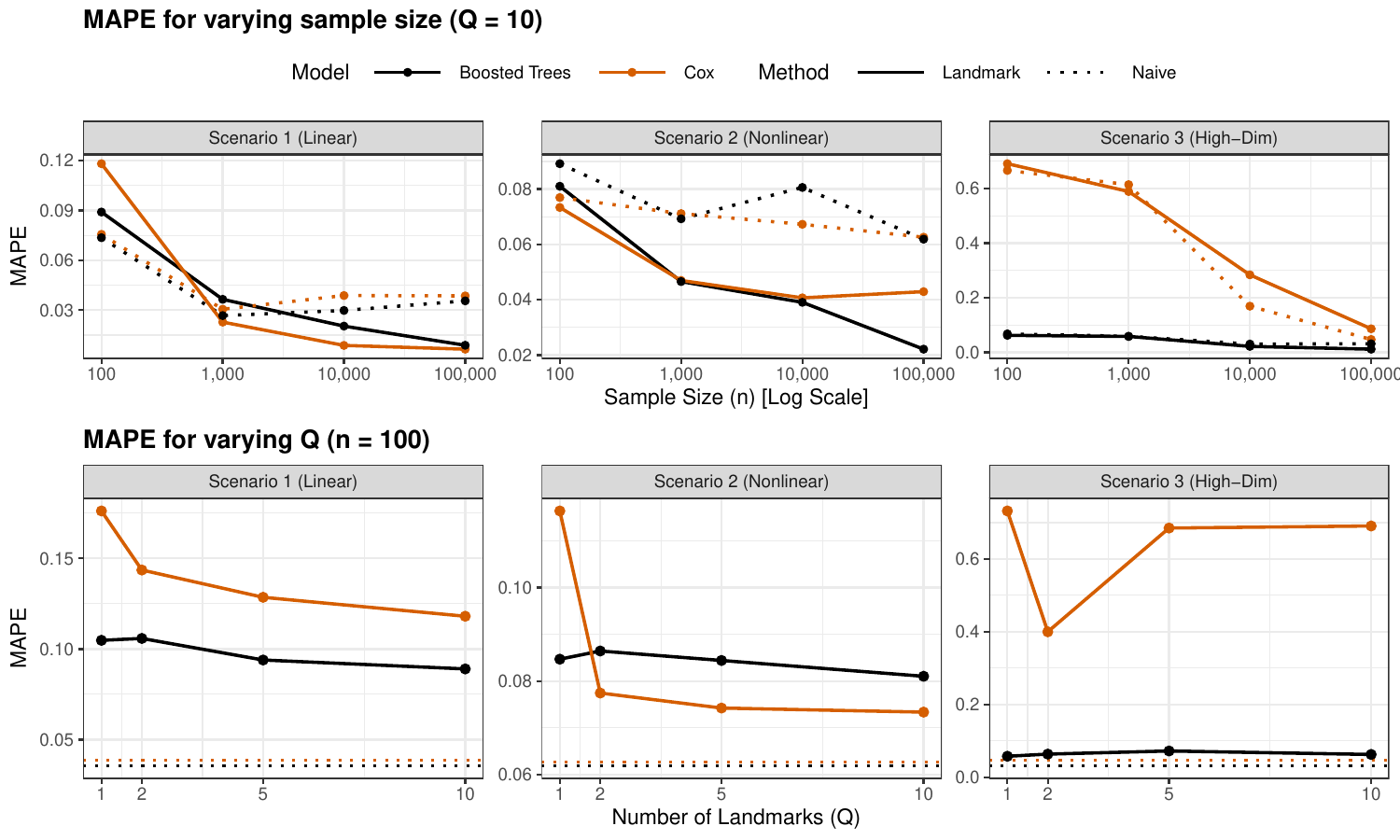}
    \caption{Simulation results. Mean absolute percentage error for the estimated future conditional survival probabilities using the proposed boosted tree landmarking supermodel (black solid lines), the Cox landmarking supermodel (orange solid lines), the naive Cox estimator (orange dotted line), and the naive boosted tree estimator (black dotted line). The top row shows performance as a function of sample size $n$ (log scale) with $Q=10$ landmarks per subject; the bottom row shows performance as a function of the number of landmarks $Q$ with sample size $n=100$. In the bottom row, the naive models appear as horizontal lines as they do not landmark the data. The columns correspond to Scenario 1 (Linear Markovian), Scenario 2 (Nonlinear Non-Markovian), and Scenario 3 (High-dimensional).}
    \label{fig:MAPE}
\end{figure}

\begin{figure}[ht!]
    \centering
\includegraphics[width=1\linewidth]{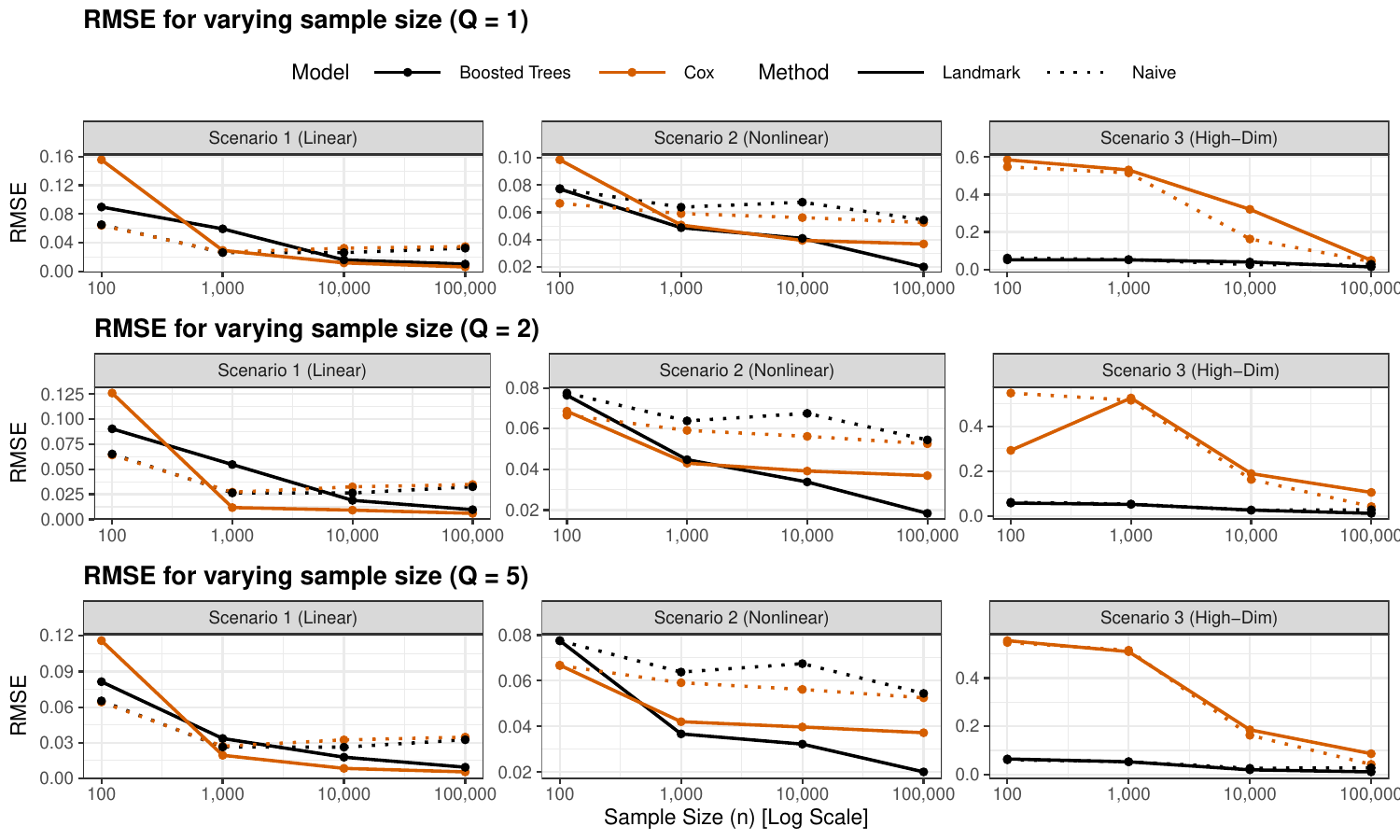}
    \caption{Simulation results. Root mean squared error for the estimated future conditional survival probabilities using the proposed boosted tree landmarking supermodel (black solid lines), the Cox landmarking supermodel (orange solid lines), the naive Cox estimator (orange dotted line), and the naive boosted tree estimator (black dotted line). The rows show performance as a function of sample size $n$ (log scale) with the top row having $Q=1$ landmarks per subject, the middle row having $Q=2$, and the bottom row having $Q=5$. The columns correspond to Scenario 1 (Linear Markovian), Scenario 2 (Nonlinear Non-Markovian), and Scenario 3 (High-dimensional).}
    \label{fig:RMSEn}
\end{figure}

\begin{figure}[ht!]
    \centering
\includegraphics[width=1\linewidth]{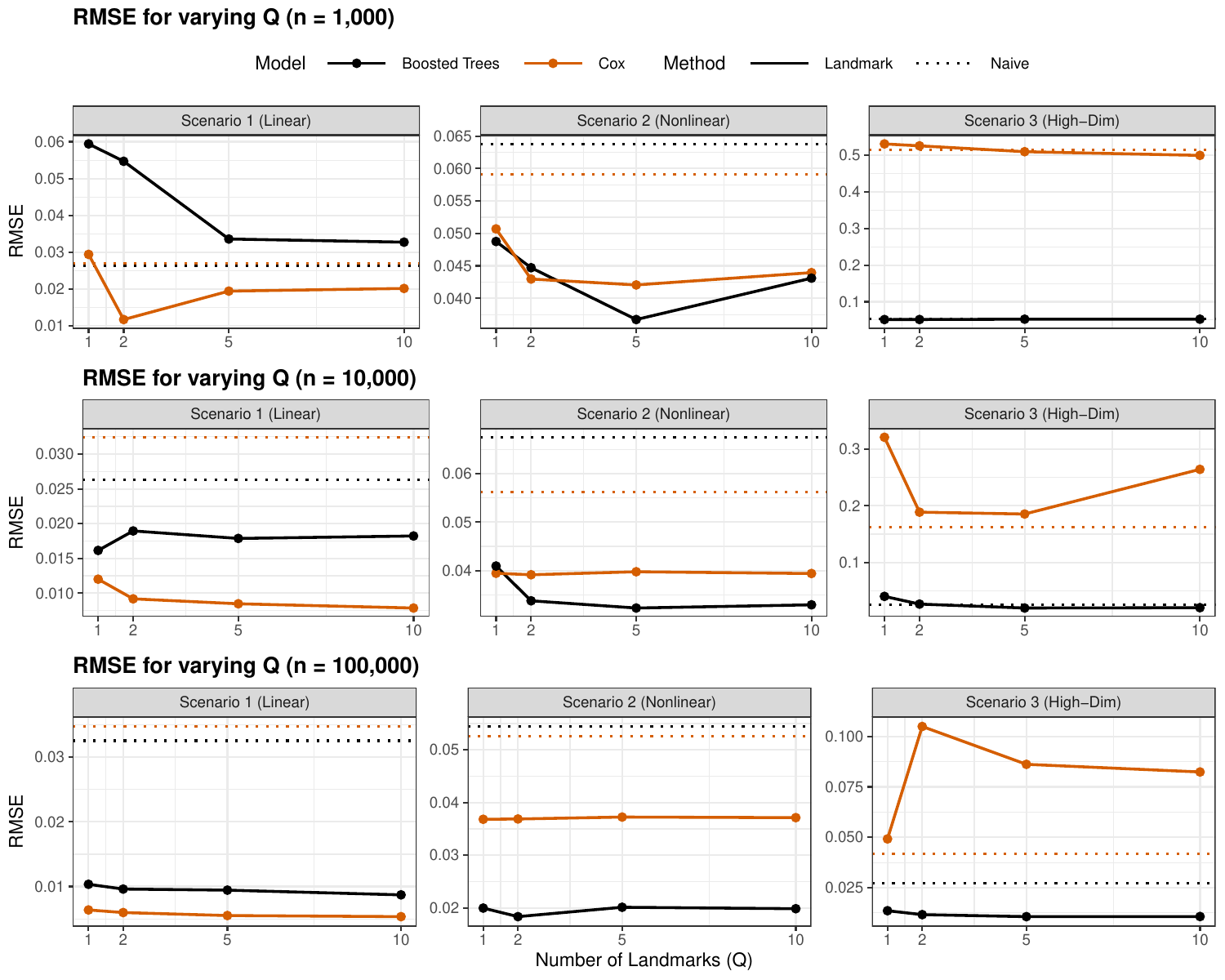}
    \caption{Simulation results. Root mean squared error for the estimated future conditional survival probabilities using the proposed boosted tree landmarking supermodel (black solid lines), the Cox landmarking supermodel (orange solid lines), the naive Cox estimator (orange dotted line), and the naive boosted tree estimator (black dotted line). The rows show performance as a function of $Q$ with the top row having sample size $n=1 \, 000$, the middle row having $n=10 \, 000$, and the bottom row having $n= 100 \, 000$. The naive models appear as horizontal lines as they do not landmark the data. The columns correspond to Scenario 1 (Linear Markovian), Scenario 2 (Nonlinear Non-Markovian), and Scenario 3 (High-dimensional).}
    \label{fig:RMSEq}
\end{figure}

\FloatBarrier

\subsection{Computational details}

Analyses were performed in R 4.2.3 on Ubuntu 22.04.5 LTS in an RStudio Server session. The machine had access to 16 logical CPU cores (host CPU: AMD EPYC 9V74 80‑Core Processor) and 126 GB RAM. In Table~\ref{tab:training_times}, we report the training times for each model. For the boosted tree models, this includes the time spent on $K$-fold cross-validation to select the number of trees. The Cox landmarking supermodel did not converge in the high-dimensional setting for $n=100$ with $Q=1$ and $Q=2$. This was the case even after the maximum number of iterations was raised to 1000. We therefore used the model that was reached after 20 iterations, which is the default maximum number of iterations in \texttt{coxph}.

\begin{table}[ht]
\centering
\def~{\hphantom{0}}
\caption{Complete training times (in seconds) for Cox and Boosted Tree based estimation.}
\label{tab:training_times}
\begin{tabular}{lrlrrrrr}
  \toprule
Scenario & $n$ & Model & Naive & $Q=1$ & $Q=2$ & $Q=5$ & $Q=10$ \\ 
  \midrule
Scenario 1 & 100 & Cox & 0.01 & 1.51 & 1.42 & 1.53 & 2.39 \\ 
   &  & Boosted Trees & 0.02 & 0.02 & 0.04 & 0.04 & 0.31 \\ 
   & 1,000 & Cox & 0.02 & 1.53 & 1.44 & 1.49 & 1.54 \\ 
   &  & Boosted Trees & 1.34 & 0.08 & 0.16 & 2.04 & 6.39 \\ 
   & 10,000 & Cox & 0.10 & 1.52 & 2.77 & 1.78 & 2.08 \\ 
   &  & Boosted Trees & 30.19 & 12.28 & 21.62 & 59.66 & 105.95 \\ 
   & 100,000 & Cox & 1.22 & 2.06 & 2.87 & 5.12 & 9.69 \\ 
   &  & Boosted Trees & 774.66 & 219.82 & 958.30 & 2398.46 & 7455.30 \\ 
   \midrule 
  Scenario 2 & 100 & Cox & 0.01 & 1.41 & 1.42 & 1.40 & 1.30 \\ 
   &  & Boosted Trees & 0.03 & 0.02 & 0.02 & 0.06 & 0.08 \\ 
   & 1,000 & Cox & 0.02 & 1.43 & 1.43 & 1.38 & 1.31 \\ 
   &  & Boosted Trees & 1.37 & 0.81 & 1.00 & 4.19 & 5.85 \\ 
   & 10,000 & Cox & 0.09 & 1.45 & 1.50 & 1.60 & 2.35 \\ 
   &  & Boosted Trees & 28.25 & 6.94 & 27.25 & 52.49 & 117.60 \\ 
   & 100,000 & Cox & 0.90 & 1.85 & 2.58 & 4.73 & 8.16 \\ 
   &  & Boosted Trees & 647.36 & 217.26 & 552.69 & 2000.24 & 4431.40 \\ 
   \midrule
  Scenario 3 & 100 & Cox & 0.02 & 1.40 & 1.48 & 1.45 & 1.40 \\ 
   &  & Boosted Trees & 0.34 & 0.34 & 0.35 & 0.39 & 0.47 \\ 
   & 1,000 & Cox & 0.07 & 1.42 & 1.40 & 1.45 & 1.59 \\ 
   &  & Boosted Trees & 0.70 & 0.51 & 0.68 & 1.24 & 2.62 \\ 
   & 10,000 & Cox & 0.47 & 1.59 & 1.86 & 2.87 & 4.65 \\ 
   &  & Boosted Trees & 265.61 & 48.56 & 103.87 & 123.31 & 215.24 \\ 
   & 100,000 & Cox & 6.35 & 4.65 & 7.56 & 18.58 & 63.05 \\ 
   &  & Boosted Trees & 1715.46 & 473.37 & 2326.20 & 4633.20 & 5950.08 \\ 
   \bottomrule
\end{tabular}
\end{table}

\section{Additional details for data application} \label{supplement:dataDetails}

\subsection{Additional figures}

\begin{figure}[ht!]
    \centering    \includegraphics[scale=0.65]{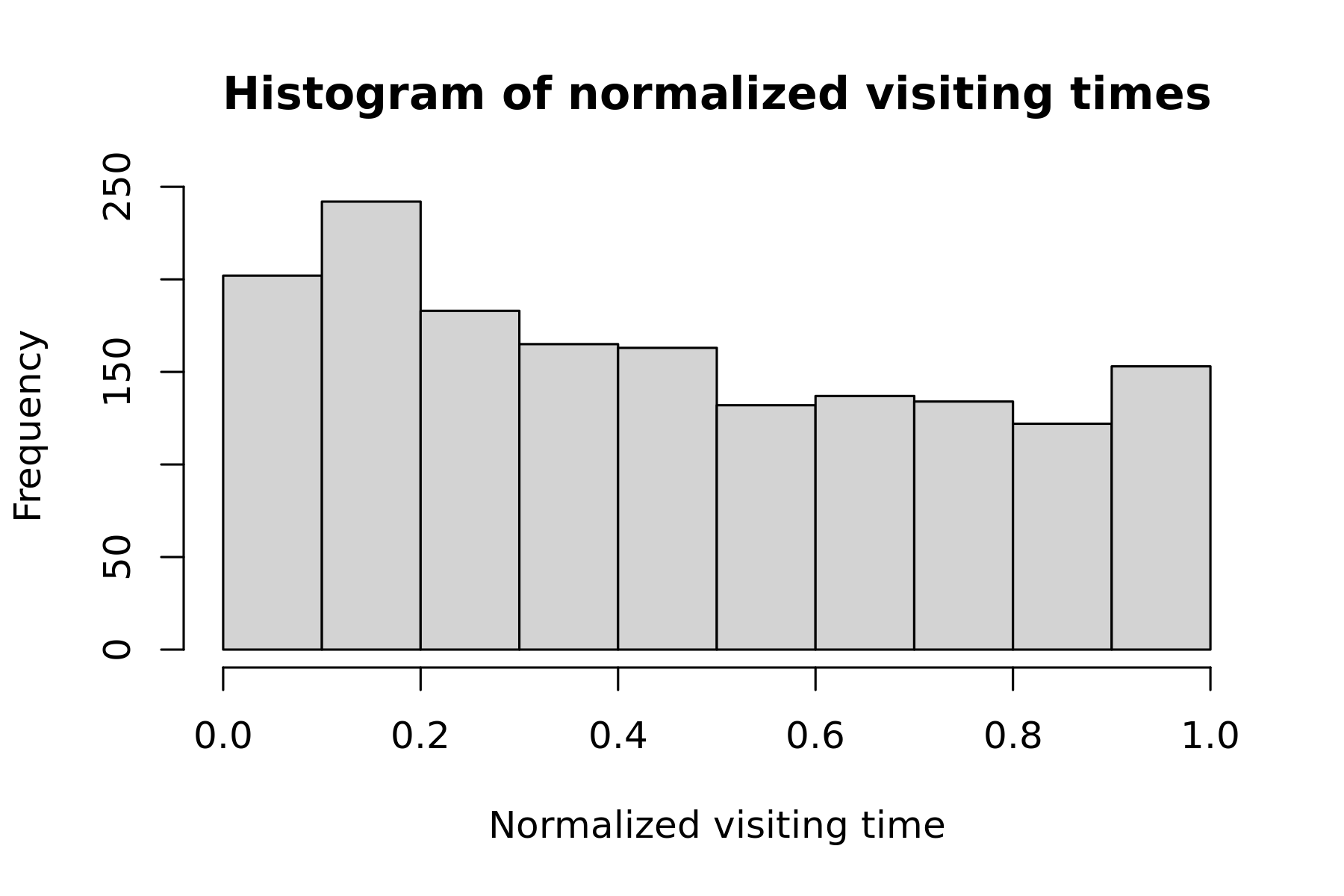}
    \caption{Histogram of visiting times normalized by the time of study exit and excluding the time of enrollment.}
    \label{fig:visitHist}
\end{figure}

Fig.~\ref{fig:visitHist} is a histogram of the post-enrollment visit times when normalized by the time of study exit. As noted in the main text, it is seen that the distribution is approximately uniform. 

\begin{figure}[ht!]
    \centering  
    \hspace{-3.7cm}
    \includegraphics[scale=0.85]{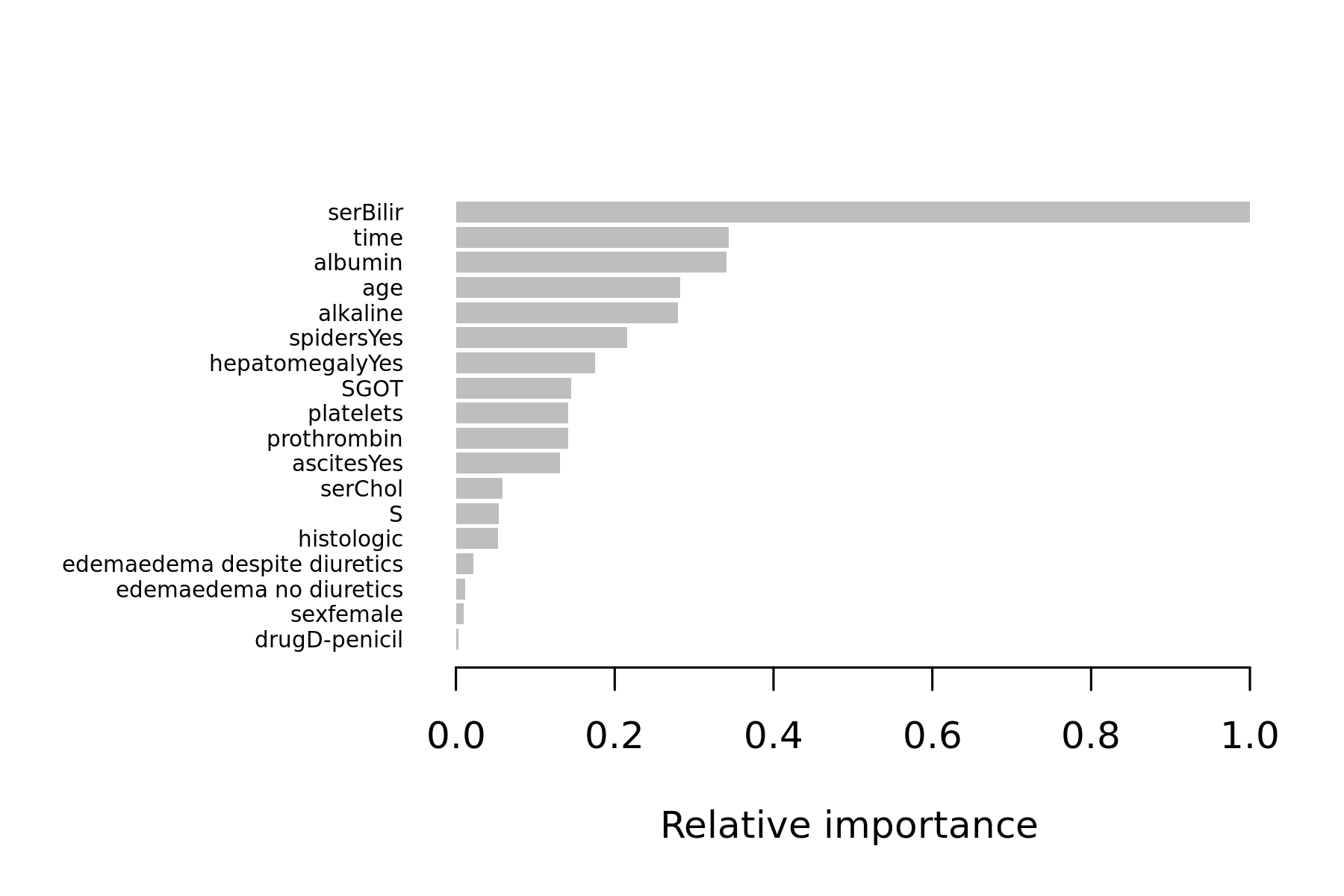}
    \caption{Normalized importance of the covariates based on the fitted model, where importance for variable $j$ is the sum of the loss reductions across every split on variable $j$. The covariates are ordered according to their importance.}
    \label{fig:importance}
\end{figure}

\begin{figure}[ht!]
    \centering    
    \hspace{-2cm}
    \includegraphics[scale=0.85]{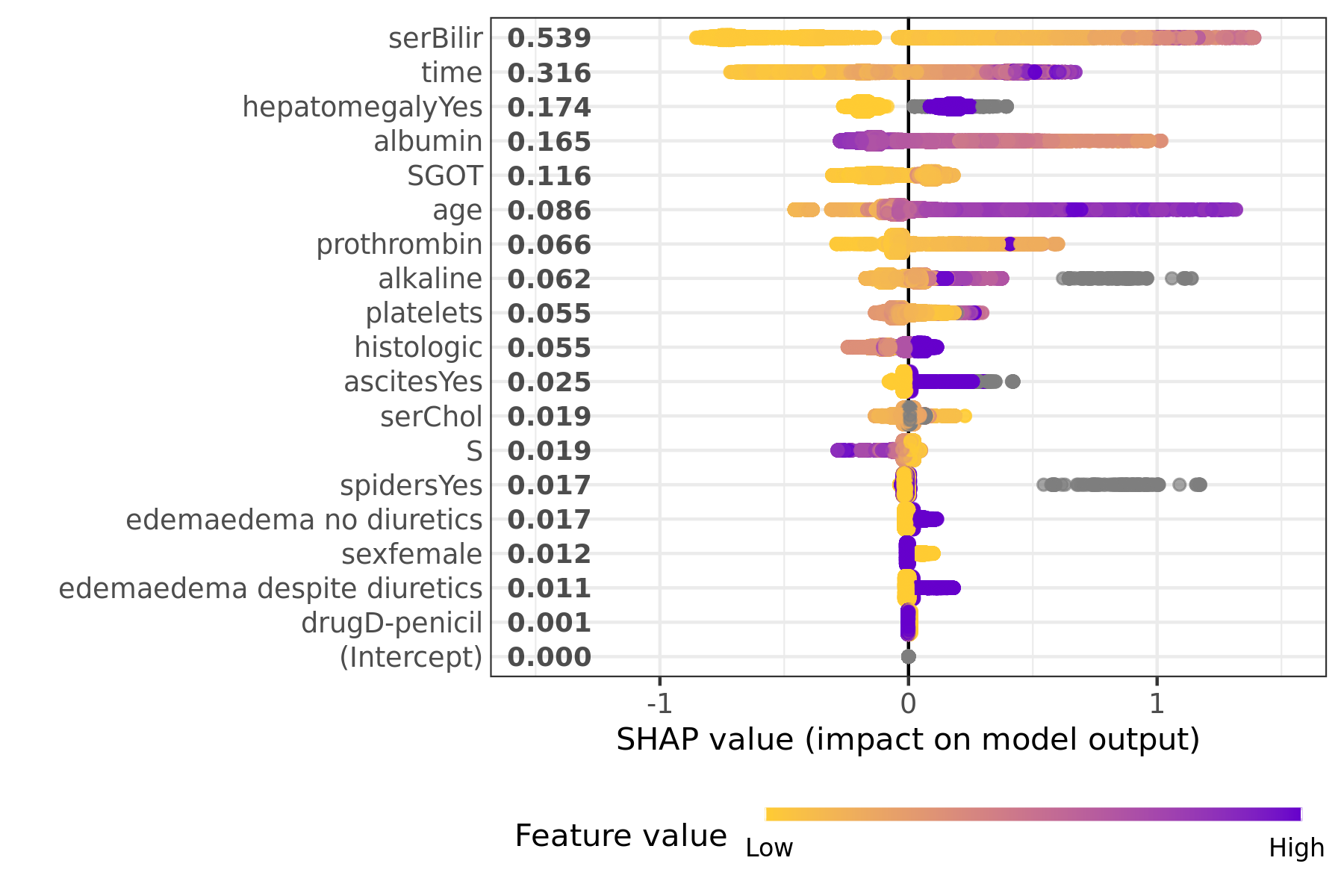}
    \caption{Shapley values (points) of the covariates for the fitted model. The number on the left indicates the mean absolute Shapley value and determines the order. The colors illustrate whether the observed value is missing (gray), low (yellow), or high (purple), where the latter two are based on the distance to the minimal and maximal observed value.}
    \label{fig:shapley}
\end{figure}

Fig.~\ref{fig:importance} is based on the native XGBoost importance measure for each covariate, which is the sum of the reduction in the loss for each split on that particular covariate. These values are normalized by the largest importance value, which in our case is the covariate serBilir that encodes the level of bilirubin. The normalized values are visualized in a barplot in descending order. 
Another widely used importance measure is based on Shapley values. Fig.~\ref{fig:shapley} depicts the Shapley values computed over all the observations in the data. The horizontal position represents the estimated contribution of the covariate to the future conditional log-hazard. Positive values indicate increased predicted risk, and negative values indicate decreased predicted risk. Colors encode the observed covariate values, ranging from low (yellow) to high (purple) based on the distance to the minimal and maximal observed value. Missing values are colored gray. The covariates are ordered by their mean absolute Shapley value, shown to the left of each row, which provides a measure of relative importance. The level of bilirubin is again identified as the most important covariate for the prediction.

Table~\ref{tab:covNA} shows the number of missing values for each covariate across the 1633 post-enrollment visit times. Covariates with no missing values are not included in the table. We use the native handling of missing values in XGBoost. An alternative would have been to exclude the covariates or exclude subjects with missing values in any of the covariates. Serveral of these are ranked relatively highly in the importance measures of Figs.~\ref{fig:importance}-\ref{fig:shapley}, and it would thus have been disadvantageous to remove them entirely.

\begin{table}[ht!]
\caption{Number of missing values across the 1633 post-enrollment visit times. Covariates not included have no missing values.}
\label{tab:covNA}
\centering
\begin{tabular}{lrrrrr}
\hline
ascites & hepatomegaly & spiders & serChol & alkaline & platelets \\
\hline
60 & 61 & 58 & 793 & 60 & 69 \\
\hline
\end{tabular}
\end{table}

\FloatBarrier

Table~\ref{tab:subjects} lists the covariates for two selected subjects with different risk profiles. Subject A has several clinical and biochemical factors that correlate with a poor survival prognosis in the medical literature, while subject B presents mild disease.

\begin{table}[ht!]
\centering
\caption{Covariate values for the two subjects with different risk profiles.}
\label{tab:subjects}
\begin{tabular}{lcc}
\toprule
Variable & Subject A (id 128) & Subject B (id 25) \\
\midrule
drug             & D-penicil            & placebo \\
age              & 43.65                & 46.16 \\
sex              & female               & female \\
ascites          & No                   & No \\
hepatomegaly     & Yes                  & No \\
spiders          & Yes                  & No \\
edema            & edema despite diuretics & No edema \\
serBilir         & 18.5                 & 0.5 \\
serChol          & 168                  & NA \\
albumin          & 3.37                 & 3.93 \\
alkaline         & 1928                 & 589 \\
SGOT             & 398.4                & 328.6 \\
platelets        & 115                  & 285 \\
prothrombin      & 13.1                 & 10.3 \\
histologic       & 4                    & 2 \\
landmark time             & 0.8515               & 0.5448 \\
\bottomrule
\end{tabular}
\end{table}

\FloatBarrier

\section{Proof of Proposition~\ref{prop:approx}} \label{supplement:space}

For simplicity, we take $\nu_{kj}=1$, $d_{kj} = p+2$, and $m_k \rightarrow \infty$. Assume also that the partitions $\mathcal{P}_k$ of $\mathcal{D}$ become increasingly fine as $k \rightarrow \infty$ in the sense that the maximal diameter of the hypercubes $B_{k\ell}$ in $\mathcal{P}_k$ tends to zero. Let 
    \[
    \tilde{\mathcal{F}} = \Bigl\{ F  : 
    F \text{ is $\mu$-measurable and } 
    \log \Lambda_L \leq F \leq \log \Lambda_U \;\mu\text{-almost surely} \Bigr\}
    \]
be the target space and $\mathcal{U}=\cup_{k=1}^\infty \mathcal{F}_k$. Recall that $\mathcal{F}$ is defined as the completion of $\mathcal{U}$ under the metric $d$.  Completions are unique up to isometries, but since we consider $\mathcal{U}$ as a subset of $L^1(\mu)$, the completion is unique and equal to the closure of $\mathcal{U}$ within $L^1(\mu)$. We prove that $\mathcal{F}=\tilde{\mathcal{F}}$ by showing that $\tilde{\mathcal{F}}$ is a complete
metric space under the metric $d$ and that $\mathcal{U}$ is a dense subset of it. This implies that the completion of $\mathcal{U}$ under $d$ is $\tilde{\mathcal{F}}$.

Since $\mu$ is a finite measure we see that $\tilde{\mathcal{F}} \subset L^1(\mu)$. The space $L^1(\mu)$ is a complete normed
vector space per the Riesz-Fischer theorem (see Theorem 6.6 in~\citet{Folland:1999}). A subset of a complete metric space is complete if and only if it is closed. To show that $\tilde{\mathcal{F}}$ is a closed subset of $L^1(\mu)$, let $(F_m)_{m \geq 1}$ be a sequence in $\tilde{\mathcal{F}}$ that converges in $L^1(\mu)$ to a limit $F \in L^1(\mu)$. By Corollary 2.32 in~\citet{Folland:1999}, there exists a subsequence $(F_{m_j})_{j \geq 1}$ with $m_1<m_2<\dots$ that converges to $F$ pointwise $\mu$-almost surely. Since every element in the subsequence satisfies the bounds $\log \Lambda_L \leq F_{m_j} \leq \log \Lambda_U$, their limit $F$ also satisfies these bounds $\mu$-almost surely. Thus, $F \in\tilde{\mathcal{F}}$ so $\tilde{\mathcal{F}}$ is closed and therefore complete.

We now show that $\mathcal{U}$ is a dense subset of $\tilde{\mathcal{F}}$, in other words, for any $F \in \tilde{\mathcal{F}}$ and $\varepsilon > 0$, there exists a function $F_k \in \mathcal{U}$ such that $d(F,F_k) < \varepsilon$. With the above choice of hyperparameters, the space $\mathcal{F}_k$ corresponds exactly to the space of simple functions with sets belonging to the collection $\mathcal{P}_k$ and coefficients belonging to the interval $[\log \Lambda_L,\log \Lambda_U]$. The set of simple functions is dense in $L^1(\mu)$ (see Proposition 6.7 in~\citet{Folland:1999}), so there exists a simple function $S = \sum_{j=1}^J a_j 1\{ \cdot \in D_j\}$ with $d(F,S) < \varepsilon/2$ where $D_j$ are sets in Borel $\sigma$-algebra on $\mathcal{D}$. Furthermore, since $\log \Lambda_L \leq F \leq \log \Lambda_U$, we may assume $a_j \in [\log \Lambda_L, \log \Lambda_U]$. If we for $k$ large enough can find a simple function $F_k \in \mathcal{F}_k \subseteq \mathcal{U}$ with $d(S,F_k) < \varepsilon/2$ then $d(F,F_k)\leq d(F,S)+d(S,F_k)<\varepsilon$ as required.

We start by showing that we can approximate individual indicator functions. For a given measurable set $A$ and $\delta > 0$, we show that there exists a $k$ and $A_k$, which is a finite union of disjoint hypercubes from the partition $\mathcal{P}_k$ such that $d(1\{\cdot \in A\}, 1\{\cdot \in A_k\}) = \mu(A \Delta A_k) < \delta$. Observe first that $\mu$ is a finite Borel measure on a subset of $\mathbb{R}^{p+2}$, so it is regular. This implies that there exists a compact set $K$ with $K \subseteq A$ and $\mu(A \setminus K) < \delta/2$ as well as an open set $O$ with $A \subseteq O$ and $\mu(O \setminus A) < \delta/2$. Since $K$ is compact and $O$ is open with $K \subset  O$, the distance between $K$ and the complement of $O$ is strictly positive. The maximal diameter of the hypercubes $B_{k\ell}$ tends to zero, so we can choose $k$ large enough  so that any  $B_{k\ell}$ that intersects $K$ must be entirely contained within $O$. Define $A_k = \cup_{\ell : B_{k\ell}  \cap K \neq \emptyset } B_{k\ell}$ and note that the union is finite since $\mathcal{P}_k$ is finite. By construction, we have $K \subseteq A_k \subseteq O$. Now consider the symmetric difference $A \Delta A_k =(A\setminus A_k) \cup (A_k \setminus A)$. Note $A_k \setminus A \subseteq O \setminus A$ so $\mu(A_k \setminus A) \leq \mu(O \setminus A) < \delta/2$. Likewise $A \setminus A_k \subseteq A \setminus K$ so $\mu(A \setminus A_k) \leq \mu(A \setminus K) < \delta/2$. This implies $\mu(A \Delta A_k) < \delta$. Note that $1\{\cdot \in A_k\} = \sum_{\ell : B_{k\ell}  \cap K \neq \emptyset} 1\{ \cdot \in B_{k\ell} \}$ and that the right-hand side can be obtained as a tree in $\mathcal{F}_k$. For $k$ large enough, we can hence find $D_{kj}$ such that $\mu(D_j \Delta D_{kj})$ is as small as we want.

We define $F_k \in \mathcal{F}_k$ as being the simple function with sets $D_{kj}$ and leaf values $a_j$. Choose $k$ large enough that $\mu(D_j \Delta D_{kj}) < \varepsilon/(2 \sum_{j=1}^J \vert a_j \vert)$. This implies $d(S,F_k) = \norm{S-F_k}_{\mu,1} \leq \sum_{j=1}^J \vert a_j \vert \times \norm{1\{\cdot \in D_j\}-1\{\cdot \in D_{kj}\}}_{\mu,1} < \varepsilon/2$. 

\section{Conditions for Assumption~\ref{assumption:maximize}} \label{supplement:maximize}

Boosting is usually formulated for real-valued outcomes rather than outcomes of stochastic processes, so in order to use results from the literature, we use the Poisson representation of the criterion function described in Section~\ref{subsec:estimation}. Associate with each set $B_{n\ell}$ a single representative value $(t_{n\ell},s_{n\ell},w_{n\ell}) \in B_{n \ell}$ and treat this as the outcome of a  covariate denoted $X_{in\ell q}$ associated with $O_{in\ell}(S_{iq},T]$. For simplicity, re-index $O_{in\ell}(S_{iq},T]$, $E_{in\ell}(S_{iq},T]$, and $X_{in\ell q}$ for varying $(i,\ell,q)$ and fixed $n$ in terms of a single index $k$ leading to $O_k$, $E_k$, and $X_k$ for $k=1,\dots,K$ where we only keep the observations where the exposure is non-zero. Note that the $n$-dependence suppressed in the notation. We may now represent the (rescaled) empirical loss-function $R_n(F) = -(nQ/K) \times M_n(F)$ as
$$R_n(F) = K^{-1} \sum_{k=1}^K \psi_k \{F(X_k),O_k\}$$
for the individual loss $\psi_k \{F(X_k),O_k\}=\exp\{F(X_k)\}E_k - F(X_k) O_k $ with $F\in \mathcal{F}_n$. Showing Assumption~\ref{assumption:maximize} holds thus corresponds to showing 
\begin{align} \label{eq:Rn}
    R_n(\hat{F}_n) \leq \inf_{F \in \mathcal{F}_n } R_n(F) + o_P(1).
\end{align}
For our boosting estimator, we use Algorithm~\ref{alg:boost} for a fixed value of $n$. The function $\textnormal{cap}$ in the algorithm is the truncation operator defined as $\textnormal{cap}(F)=\max\{\min\{F,\log \Lambda_U\},\log \Lambda_L\}$ and $$\xi_k(x,y)=\frac{\partial}{\partial x} \psi_k(x,y) = \exp(x)E_k-y$$ 
is the gradient. 

 \begin{algorithm}[ht!] 
\caption{Gradient Boosting}
\label{alg:boost}
\KwIn{Step size $\nu>0$.}
\KwData{Initial function $F_{0} \in \mathcal{T}(d,\mathcal{P}_n)$ and iteration counter $j \leftarrow 0$.}
\Repeat{
  Compute the tree that best approximates the gradient in squared error:
  \[
    g_{j+1} \in \argmin_{g \in \mathcal{T}(d,\mathcal{P}_n)} K^{-1} \sum_{k=1}^K [-\xi_k\{F_j(X_k),O_k\}-g(X_k)]^2.
  \]
  Update the predictor:
  \[
     F_{j+1} \leftarrow F_j + \nu g_{j+1}.
  \]
  Set $j \leftarrow j + 1$\;
}
Set $\hat{F}_n \leftarrow  \textnormal{cap}(F_{\infty})$.
\end{algorithm}

Except for the final truncation, the boosting algorithm is identical to Algorithm 2 in~\citet{Biau.Cadre:2021} in the case where the gradient exists and the loss function is the empirical loss function. Algorithm 2 in~\citet{Biau.Cadre:2021} corresponds to Algorithm 1 in~\citet{Friedman:2001} except that a fixed stepsize is used. 

We proceed by adapting the results of~\citet{Biau.Cadre:2021} to our setting in order to use their Theorem 3.2, which states that the boosting estimator minimizes the loss. We accommodate (1) that the space for the individual trees has to be the same across boosting iterations and be closed under multiplication by scalars, (2) that the individual loss $\psi_k$ can be negative, (3) that $\psi_k$ are different for different $k$ due to the influence of exposure, and (4) that the final boosting estimator is truncated.

(1): Set $(m_k,d_{kj},\nu_{kj}) = (\infty,d,\nu)$ for fixed $d \in \mathbb{N}$ and $\nu > 0$. The space for the individual trees is thus $\mathcal{T}(d,\mathcal{P}_n)$. As required by the setup of Algorithm 2 in~\citet{Biau.Cadre:2021} the space is the same for each boosting iteration, the space is closed under multiplication by scalars, the stepsize is fixed, and the boosting algorithm runs indefinitely. The function space that the boosting estimator belongs to prior to truncation is thus
$$\Tilde{\mathcal{F}}_n=\left\{ F=\sum_{j=1}^{\infty} \nu g_j : g_j \in \mathcal{T}(d,\mathcal{P}_n)  \right\},$$
which is the linear span of $\mathcal{T}(d,\mathcal{P}_n)$ and the unbounded version of $\mathcal{F}_n$.

(2): Since $F \in \mathcal{F}_n$ is bounded within $[\log \Lambda_L,\log \Lambda_U]$, the individual loss $\psi_k(F,O_k)$ is uniformly bounded below by $-O_k \log \Lambda_U$ as a function of $F$. The only place where non-negativity is used for obtaining Theorem 3.2 in~\citet{Biau.Cadre:2021} is in Lemma 3.2 to conclude that since the sequence of loss functions evaluated at the boosting iterates $F_j$ is decreasing, it must converge. This also holds in our setting, where the loss is bounded below as a function of $F$.

(3) and (4): We treat these together by modifying and verifying assumption (A1)-(A3) from~\citet{Biau.Cadre:2021}, which together imply Theorem 3.2. The individual loss $\psi_k$ are not natively $\alpha$-strongly convex since the second derivatives are $\exp(x)E_k$ which goes to $0$ as $x \rightarrow -\infty$. Similarly, the gradient $\xi_k$ is not Lipschitz continuous due to the $\exp(x)E_k$ term. To remedy this, we assume that the boosting iterates stays confined to a bounded set with high probability.

\begin{assumption} \label{assumption:boostingBound}
Let 
$$S_n = (\forall j: F_j \in [\log \Lambda_L,\log \Lambda_U])$$
be the event that the boosting iterates $F_j$ all take values within the interval $[\log \Lambda_L,\log \Lambda_U]$. We assume that $P(S_n) \rightarrow 1$ for $n \rightarrow \infty$. Furthermore, we assume $R_n(\hat{F}_n) -  \inf_{F \in \mathcal{F}_n } R_n(F) = O_P(1).$
\end{assumption}

 Furthermore, without loss of generality, we assume $0 \in [\log \Lambda_L,\log \Lambda_U]$. Assumption~\ref{assumption:boostingBound} is mild in practice since the bounds can be made as wide as needed. Under this assumption, we only need to show Equation~\eqref{eq:Rn} on the event $S_n$ since 
 $$ 1\{S_n^c\} \times \{R_n(\hat{F}_n) -  \inf_{F \in \mathcal{F}_n } R_n(F) \} \leq 1\{S_n^c\} \times \vert R_n(\hat{F}_n) -  \inf_{F \in \mathcal{F}_n } R_n(F) \vert = o_p(1)$$
 since $\vert R_n(\hat{F}_n) -  \inf_{F \in \mathcal{F}_n } R_n(F) \vert$ is tight and $P(S_n^c) \rightarrow 0$ for $n \rightarrow \infty$.

For the empirical setting of ~\citet{Biau.Cadre:2021}, one may check that if the boosting iterates all take values within some bounded interval, then Theorem 3.2 follows if Assumption (A1)-(A3) hold for $x$ in this bounded interval rather than for all of $\mathbb{R}$. On $S_n$, we thus only need to consider $\psi_k(x,y)$ for $x \in [\log \Lambda_L,\log \Lambda_U]$. Furthermore, the relevant proofs all hold when $\psi_k$  differ across $k$ if Assumption (A1)-(A3) hold uniformly over $k$. Since the truncation in Algorithm~\ref{alg:boost} is superfluous on $S_n$, Theorem 3.2 of~\citet{Biau.Cadre:2021} holds if we show these slightly stronger conditions:
\begin{itemize}
    \setlength{\itemindent}{2em}
    \item[(A1)] $\psi_k(0,O_k) < \infty$ and $\sup_{x \in [\log \Lambda_L,\log \Lambda_U]}\xi_k(x,O_k) < \infty$ for all $k$.
    \item[(A2)] $\frac{\partial}{\partial x} \xi_k(x,y) \geq  \alpha >0$ for all $k$ and $x \in [\log \Lambda_L,\log \Lambda_U]$.
    \item[(A3)] $\vert \xi_k(x_1,y) - \xi_k(x_2,y) \vert \leq \delta \vert x_1-x_2 \vert$ for a $\delta >0$ and all $k,y,$ and $(x_1,x_2) \in [\log \Lambda_L,\log \Lambda_U]^2$.
\end{itemize}
Condition (A1) holds trivially. Condition (A2) holds since only observations with strictly positive exposures are included, and $\Lambda_L > 0$. For condition (A3), note 
\begin{align*}
    \vert \xi_k(x_1,y) - \xi_k(x_2,y) \vert \leq E_k \Lambda_U \vert x_1-x_2 \vert 
\end{align*}
using that $\vert \exp(x_1)-\exp(x_2) \vert \leq \exp(\max\{x_1,x_2\})\vert x_1 - x_2 \vert$ as shown under Condition 3.1 (ii) in Section E and that $(x_1,x_2) \in [\log \Lambda_L,\log \Lambda_U]^2$. We may hence use $\delta = \inf_k E_k \Lambda_U > 0$. 

By Theorem 3.2 in~\citet{Biau.Cadre:2021}, we thus get 
$$1\{S_n\} \times R_n(\hat{F}_n) = 1\{S_n\} \times \inf_{F \in \tilde{\mathcal{F}}_n} R_n(F) \leq 1\{S_n\} \times \inf_{F \in \mathcal{F}_n} R_n(F)$$
which implies that Assumption~\ref{assumption:maximize} holds.

\section{Proof of Theorem~\ref{theorem:consistency}} \label{supplement:consistency}

The strategy is to verify the conditions of Theorem 3.1 in~\citet{Chen:2007}. Many of these conditions are statements about the population criterion function, so we begin by deriving a more convenient representation. Note that 
$$E\left[ \int_0^T \int_s^T F\{t,s,W(s)\} \diff N(t) \diff s \right ] = E\left[ \int_0^T \int_s^T F\{t,s,W(s)\} Y(t) \lambda\{t,s,W(s)\} \diff t \diff s \right ]$$
since $t \mapsto \lambda\{t,s,W(s)\}Y(t)$ is the intensity process for $N$ under the filtration $t \mapsto \mathcal{G}_{t,s}$ and using that the integrand $F$ is predictable with respect to this filtration. We thus get
$$M(F)=\int F \lambda -\exp(F) \diff \mu$$
This way of representing the population likelihood loss is very similar to the one in~\citet{Lee.etal:2021}, but the reference measures are different. By assuming that the function space is a linear span of a finite number of basis functions (which is true for $\mathcal{F}_k$ but not necessarily for $\mathcal{F}$) they also get a convenient representation for the empirical risk, which does not seem to have a natural counterpart in our setting as we do not make this restriction on the function space.

We define the domain of $M$ to be the space of bounded measurable functions from $\mathcal{D}$ to $\mathbb{R}$. This only has implications for the calculations leading to Condition 3.1 (ii), as all the other arguments only use $M$ evalued in functions belonging to $\mathcal{F}$. 

\begin{remark}
      Although not directly useful as it follows from Condition 3.1 (ii) below, we provide a simple proof that $M$ is maximized at $\log \lambda$.  It suffices to prove that $F\lambda - \exp(F) \leq \log \lambda \times \lambda - \lambda$ for any $F \in \mathbb{R}$ and $\lambda > 0$. Observe $\frac{\diff}{\diff F}\{F\lambda - \exp(F) \} =\lambda-\exp(F)$, so the first order condition implies $F=\log \lambda$ and it is a maximum since $\frac{\diff^2}{\diff F}\{F\lambda - \exp(F) \} =-\exp(F)<0$. Thus, by bounding the integrand pointwise, we obtain $M(F) \leq M(\log \lambda)$. Since $\log \lambda \in \mathcal{F}$, it is the maximizer of $M$ over $\mathcal{F}$.
\end{remark}
 
We now verify Condition 3.1-3.5 of~\citet{Chen:2007}:

\noindent \textbf{Condition 3.1 (i)}: Recall that we have assumed $\log \lambda  \in \mathcal{F}$. Therefore, this condition holds by the uniform boundedness imposed in Assumption~\ref{assumption:Bounded} combined with the observation that $\mu$ is a finite measure since it is trivially bounded above by $T$.

 \noindent \textbf{Condition 3.1 (ii)}: Note that 
\begin{align*}
    M(\log \lambda )-\sup_{F\in\mathcal{F}_k : d(F,\log \lambda ) \geq \varepsilon} M(F) &\geq M(\log \lambda )-\sup_{F\in\mathcal{F} : d(F,\log \lambda ) \geq \varepsilon} M(F)\\
    &= \inf_{F\in\mathcal{F} : d(F,\log \lambda ) \geq \varepsilon} \{M(\log \lambda ) -M(F)\}.
\end{align*}
By an interchange argument, we have for bounded measurable functions $F$ and $f$ that
$$\frac{\diff}{\diff \theta}M(F+\theta f) = \int f\lambda -\exp(F+\theta  f)  f \diff \mu,$$
$$\frac{\diff^2}{\diff \theta^2}M(F+\theta f) = -\int \exp(F+\theta  f)  f^2 \diff \mu.$$
 Let $G(\theta)=M\{\log \lambda +\theta  (F-\log \lambda )\}$ and note $G(0)=M(\log \lambda )$ and $G(1)=M(F)$. With a Taylor expansion of $G$ around $\theta=0$ evaluated at $\theta=1$ we obtain
\begin{align*}
    G(1) - G(0) &= G'(0)  (1-0) + \frac{1}{2}G''(p_F)(1-0)^2 \\
    &= \int (F-\log \lambda )\exp(\log \lambda )-\exp(\log \lambda )(F-\log \lambda ) \diff \mu \\
    &\quad-\frac{1}{2}\int \exp\{\log \lambda +p_F(F-\log \lambda )\}(F-\log \lambda )^2 \diff \mu \\
    &=  -\frac{1}{2} \int \exp(\log \lambda )^{1-p_F}\exp(F)^{p_F}(F-\log \lambda )^2 \diff \mu 
\end{align*}
for some $p_F \in (0,1)$. This implies
\begin{align*}
     M(\log \lambda )-M(F) &= -\{G(1)-G(0)\} \\
     &\geq \frac{1}{2}\Lambda_L \int (F-\log \lambda )^2 \diff \mu.
\end{align*}
By Cauchy-Schwartz, we see
\begin{align*}
    \left(\int \vert F-\log \lambda  \vert \times 1 \diff \mu \right)^2  \leq \int (F-\log \lambda )^2 \diff \mu \times \int 1^2 \diff \mu \leq T \int (F-\log \lambda )^2 \diff \mu
\end{align*}
Then we obtain
\begin{align*}
    \inf_{F\in\mathcal{F} : d(F,\log \lambda ) \geq \varepsilon} \{M(\log \lambda ) -M(F)\} \geq (2T)^{-1}\Lambda_L \varepsilon^2 > 0
\end{align*}
so we can choose $\delta(k)=1$ and $g(\varepsilon)=(2T)^{-1}\Lambda_L \varepsilon^2$.

\noindent\textbf{Condition 3.2}: The condition $\mathcal{F}_k \subseteq \mathcal{F}_{k+1} \subseteq \mathcal{F}$ holds by construction. For any $F \in \mathcal{F}$, we now construct a sequence $(F_k)_{k \geq 1}$ with $F_k \in \mathcal{F}_k$ such that $d(F_k,F) \rightarrow 0$ for $k \rightarrow \infty$. 

By Theorem 1.6.2 of~\citet{Kreyszig:1991}, we have that $\cup_{k=1}^\infty \mathcal{F}_k$ is a dense subset of $\mathcal{F}$. Therefore, there exists $(g_m)_{m \geq 1}$ so that $g_m \in \cup_{k=1}^\infty \mathcal{F}_k$ and $d(g_m, F) < m^{-1}$ for all $m$. By definition of the union, for each $m$ there exists an integer $N_m$ such that $g_m \in \mathcal{F}_{N_m}$. We can refine this sequence of indices to be strictly increasing by defining a new sequence $(\tilde{N}_m)_{m \geq 1}$: set $\tilde{N}_1 = N_1$ and let $\tilde{N}_{m+1} = \max(N_{m+1}, \tilde{N}_m + 1)$ for $m \geq 1$. This construction ensures that the sequence $(\tilde{N}_m)_{m \geq 1}$ is strictly increasing and that $g_m \in \mathcal{F}_{\tilde{N}_m}$ for all $m$ since $\mathcal{F}_{N_m} \subseteq \mathcal{F}_{\tilde{N}_m}$.

For $k < \tilde{N}_1$, we take $F_k$ to be an arbitrary element of $\mathcal{F}_k$. For $k \geq \tilde{N}_1$, we define $F_k = g_{m(k)}$, where $m(k) = \max\{m \in \mathbb{N} : \tilde{N}_m \leq k\}$. Note that the set in the definition of $m(k)$ is non-empty since $1$ is an element and it is finite since $\tilde{N}_m$ is strictly increasing, so the expression is well-defined. Note that $F_k=g_{m(k)} \in \mathcal{F}_{\tilde{N}_{m(k)}}\subseteq \mathcal{F}_k$. It remains to be shown that $d(F_k, F) \rightarrow 0$.  Observe that for any integer $M$, it holds for $k \geq \tilde{N}_M$ that $M \in \{m \in \mathbb{N} : \tilde{N}_m \leq k\}$ which implies $m(k) \geq M$. Since $M$ was arbitrary and since $m(k)$ is non-decreasing, this shows that $m(k) \rightarrow \infty$ as $k \rightarrow \infty$. Thus, $d(F_k, F) = d(g_{m(k)}, F) \rightarrow 0$.

\noindent\textbf{Condition 3.3 (i)}: Continuity implies upper semicontinuity, and continuity on $\mathcal{F}$ implies continuity on $\mathcal{F}_k$. Since $(\mathcal{F},d)$ is a metric space, sequential continuity is equivalent to continuity. Take $(F_m)_{m \geq 1}$ with $F_m \rightarrow F$ in $(\mathcal{F},d)$. Note that
$$M(F)-M(F_m) = \int (F-F_m) \lambda -\{\exp(F)-\exp(F_m)\} \diff \mu.$$
By the triangle inequality
\begin{align*}
    \vert M(F)-M(F_m)\vert &\leq \int \vert F-F_m \vert \lambda + \vert \exp(F)-\exp(F_m) \vert \diff \mu.
\end{align*}
It holds that $\vert \exp(x)-\exp(y) \vert \leq \exp\{ \max(x,y)\}  \times \vert x-y \vert$. To show this, assume WLOG that $x>y$. Then \begin{align*}
    \vert \exp(x)-\exp(y) \vert &= \exp(x)\{1-\exp(y-x)\} \\
    &= \exp(x)\left[1-\left\{1+(y-x)+\frac{z_x^2}{2} \right\} \right] \\
    &\leq \exp(x)(x-y) \\
    &= \exp\{ \max(x,y)\}\vert x-y \vert
\end{align*}
for $z_x \in (y-x,0)$ by Taylor's theorem. We thus get
\begin{align*}
    \vert M(F)-M(F_m)\vert &\leq \Lambda_U d(F,F_m).
\end{align*}
Since $d(F,F_m) \rightarrow 0$ we get $M(F_m) \rightarrow M(F)$ which implies that $M$ is continuous on $(\mathcal{F},d)$.

\noindent\textbf{Condition 3.3 (ii)}:  For $\log \lambda  \in \mathcal{F}$, choose a sequence $(F_k)_{k \geq 1}$ with $F_k \in \mathcal{F}_k$ and $d(\log \lambda ,F_k) \rightarrow 0$ which exists per Condition 3.2. We have to show $\vert M(\log \lambda )-M( F_k) \vert = o(1)$. This is an immediate consequence of the proof of Condition 3.3 (i) which showed that $M$ is continuous on $(\mathcal{F},d)$.

\noindent\textbf{Condition 3.4}: Let $q_k = \vert \mathcal{P}_k \vert$ be the cardinality of $\mathcal{P}_k$ which is finite for all $k$. We proceed by showing that $\mathcal{F}_k$ equipped with the metric $d$ is the image of a compact subset of $\mathbb{R}^{q_k}$ under a continuous function. This implies that $\mathcal{F}_k$ is compact.

Any $F \in \mathcal{F}_k$ can be uniquely represented by its coefficient vector $c_k = (c_{k\ell})_\ell \in \mathbb{R}^{q_k}$. We show that the set of permissible coefficients $Q_k$ is compact by showing it is closed and bounded and applying the Heine-Borel theorem.

Since the sets in $\mathcal{P}_k$ are disjoint, Assumption~\ref{assumption:Bounded} implies that $Q_k \subseteq [\Lambda_L,\Lambda_U]^{q_k}$ so $Q_k$ is bounded. For closedness of $Q_k$, note that it is the intersection of two sets, one related to structural constraints and the other to bounding constraints. For the structural constraints, note that by varying the coefficients for a given tree structure (i.e., a fixed partition of $\mathcal{D}$ using sets from $\mathcal{P}_k$), one obtains a linear subspace of $\mathbb{R}^{q_k}$. Because $q_k$ is finite and the depths $d_{kj}$ for varying $j$ are finite, the total number of possible tree structures in $\mathcal{F}_k$ is finite. The set of all permissible coefficient vectors under the structural constraints can therefore be written as a finite union of linear subspaces, where the union is taken over all combinations of tree structures for the $m_k$ trees. This is closed since closedness is preserved under finite unions. The bounding constraints state that $\Lambda_L \leq c_{k\ell} \leq \Lambda_U$ for any $c_k \in Q_k$. This defines a closed polyhedron in $\mathbb{R}^{q_k}$. Since $Q_k$ is the intersection of two closed sets, it is closed.

We now define a map which takes an element in $Q_k$ to its corresponding function in $\mathcal{F}_k$. It is easy to see that this function is Lipschitz continuous with respect to the uniform norm with Lipschitz constant $1$. Since all norms on finite-dimensional spaces are equivalent, the map is also Lipschitz continuous with respect to the $L^1(\mu)$ norm.

\noindent\textbf{Condition 3.5 (i)}: By Section 3.2 in Van der Vaart and Wellner, this holds if and only if the class $\{m(F) : F \in \mathcal{F}_k\}$ is Glivenko-Cantelli. We proceed via example 19.8 in Van der Vaart. Let $d_{\infty}(F_1,F_2)=\norm{F_1-F_2}_{\infty}$ be the metric induced by the supremum norm. The proof of Condition 3.4 also implies that $\mathcal{F}_k$ is compact under $d_\infty$. 

It remains to be shown the class has an integrable envelope and $F \mapsto m(F)$ is continuous surely under $d_\infty$ for $F \in \mathcal{F}_k$. 

For integrability, note by repeated use of the triangle inequality
$$\vert m(F) \vert \leq \max(\vert \log\Lambda_L\vert,\vert \log\Lambda_U\vert) N(T)+T\Lambda_U.$$
The right hand side does not depend on $F$ and is integrable, so it can be used as an integrable envelope. 

Note that
\begin{align*}
    \vert m(F)-m(F_m) \vert &\leq \int_0^T \int_s^T \vert F\{t,s,W(s)\}-F_m\{t,s,W(s)\} \vert N(\diff t) \\
    &\quad + \int_s^T Y(t) \vert \exp(F)-\exp(F_m) \vert \diff t  \diff s \\
    & \leq (T+T^2\Lambda_U)\norm{F-F_m}_{\infty}
\end{align*}
where we in the final inequality used that $N$ is a one-jump process as well as the inequality for differences of exponential functions derived under Condition 3.3 (i). Since the right hand side does not depend on the realization of the data, we have that $m$ is surely continuous under $d_\infty$.

\noindent\textbf{Condition 3.5 (ii)}: Since $\delta(k)=1$ for all $k$, this condition is equivalent to $\sup_{F \in \mathcal{F}_k} \mid M_n(F)-M(F) \mid = o_P(1)$ which is exactly Condition 3.5 (i) so it holds.

\noindent\textbf{Condition 3.5 (iii)}: This is satisfied since for any $\eta_k =o(1)$ the implication is $M_n(\hat{F}_n) \geq \sup_{F \in \mathcal{F}_n} M_n(F)-o_P(1)$ which holds by Assumption~\ref{assumption:maximize}.

By Theorem 3.1 of~\citet{Chen:2007}, we get $d(\hat{F}_n,\log \lambda )=o_P(1)$. This concludes the proof.